\documentclass{article}
 
\usepackage{amsfonts}
\usepackage{amsthm}
\usepackage{amsmath}
\usepackage{amssymb}
\usepackage{latexsym}
\usepackage[all]{xy}

\newcounter{zlist}
\newenvironment{zlist}{\begin{list}{\rm(\arabic{zlist})}{
\usecounter{zlist}\leftmargin2.5em\labelwidth2em\labelsep0.5em
\topsep0.6ex\itemsep0.3ex plus0.2ex minus0.3ex
\parsep0.3ex plus0.2ex minus0.1ex}}{\end{list}}
 
\newcounter{blist}
\newenvironment{blist}{\begin{list}{\rm(\alph{blist})}{
\usecounter{blist}\leftmargin2.5em\labelwidth2em\labelsep0.5em
\topsep0.6ex \itemsep0.3ex plus0.2ex minus0.3ex
\parsep0.3ex plus0.2ex minus0.1ex}}{\end{list}}
 
\newcounter{rlist}
\newenvironment{rlist}{\begin{list}{\rm(\roman{rlist})}{
\usecounter{rlist}\leftmargin2.5em\labelwidth2em\labelsep0.5em
\topsep0.6ex\itemsep0.3ex plus0.2ex minus0.3ex
\parsep0.3ex plus0.2ex minus0.1ex}}{\end{list}}

\newcommand{\lra}{\longrightarrow}

\newcommand{\LRa}{\Leftrightarrow}

\newcommand{\cM}{{\bf M}}

\newcommand{\Hom}{{\rm Hom}}
\newcommand{\CHom}{{^\cC{\rm Hom}}}  
\newcommand{\HomC}{{\rm Hom}^\cC}
\newcommand{\wHomC}{\widehat{\rm H}{\rm om}^\cC}
\newcommand{\wHom}{\widehat{\rm H}{\rm om}}

\newcommand{\uDelta}{{\underline{\Delta}}}
\newcommand{\ueps}{{\underline{\varepsilon}}}

\newcommand{\eC}{{\ueps_\cC}}

\newcommand{\Ke}{{\rm Ke}\,}

\newcommand{\End}{{\rm End}}
\newcommand{\EndC}{{\rm End}^\cC}

\newcommand{\Tr}{{\rm Tr}}

\newcommand{\sM}{\mbox{$\sigma [M]$}}

\newcommand{\MC}{\cM^\cC}
\newcommand{\CM}{{^\cC\cM}}

\newcommand{\Ra}{\Rightarrow}

\newcommand{\wS}{{\widehat S}}

\newcommand{\cC}{{\cal C}}

\newcommand{\Z}{\mathbb{Z}}

\newcommand{\N}{\mathbb{N}}

\newcommand{\rhu}{\scr{\rightharpoonup}}

\newcommand{\ut}{\otimes} 
\newcommand{\ot}{\otimes}

\newcommand{\scr}[1]{\mbox{\scriptsize $#1$}} 
\newcommand{\vareps}{\varepsilon}

\newcommand{\id}{I}

\newcommand{\sP}{{\sigma [P]}}

\newcommand{\roM}{\varrho^M}

\newcommand{\roN}{\varrho^N}

\newcommand{\roL}{\varrho^L}

\newcommand{\otC}{\otimes^\cC}

\def\sw#1{{\sb{\underline{#1}}}}
 
\def\proof{{\bf Proof.~}}
\def\endproof{\hfill\hbox{$\sqcup$}\llap{\hbox{$\sqcap$}}\medskip}
\def\<{{\langle}}
\def\>{{\rangle}}

\def\eps{\varepsilon}

\def\beq{\begin{equation}}
\def\eeq{\end{equation}}

\def \eC{\underline{\eps}}

\newtheorem{thm}{}[section]

\def\Label{\label}

\begin{document}

\title{\bf On Galois Comodules
}  

\author{
Robert Wisbauer \\ 
  D\"usseldorf, Germany \\ 
 }  
\date{}
\maketitle

\begin{abstract} 
Generalising the notion of Galois corings, Galois comodules
were introduced as comodules $P$ over an $A$-coring $\cC$
for which $P_A$ is finitely generated and projective and the 
evaluation map $\mu_\cC:\Hom^\cC(P,\cC)\ot_SP\to \cC$ is an isomorphism
(of corings) where $S=\End^\cC(P)$. 
It was observed that for such comodules the 
functors $\Hom_A(P,-)\ot_SP$ and $-\ot_A\cC$ from the category 
of right $A$-modules to the category of right $\cC$-comodules are isomorphic. 
In this note we call modules $P$ with this property {\em Galois
comodules} without requiring $P_A$ to be finitely generated and projective.
This generalises the old notion with this name but we show that 
essential properties and relationships are maintained.  
These comodules are close to being generators and have some 
common properties with tilting (co)modules.  
Some of our results also apply 
to generalised Hopf Galois (coalgebra Galois) extensions.  

AMS Classicfication: 16W30, 16D90, 16D80
 
\footnote{Date: September 11, 2004} 
\end{abstract}

\section{Introduction}

Let $\cC$ be a coring over the ring $A$ and put $S=\End^\cC(A)$. 
A grouplike element $g\in \cC$ makes $A$ a right $\cC$-comodule
by the coaction $\varrho^A: A\to A\ot_A\cC, \; a\mapsto 1\ot ga$. 
The notion of {\em Galois corings} $(\cC,g)$ was introduced in 
Brzezi\'nski \cite{Brz:str} by requiring the canonical map,  
$$\chi:A\ot_SA\to \cC, \quad a\ot a'\mapsto aga',$$
to be an isomorphism (of corings). It was pointed out in \cite{WiGal}
that this can be seen as the evaluation map 
 $$\mu_\cC:\Hom^\cC(A_g,\cC)\ot_S A \to \cC,$$
and that it implies bijectivity of
 $$\mu_N:\Hom^\cC(A_g,N)\ot_S A \to N,$$
for every $(\cC,A)$-injective comodule $N$.
 
The notion of Galois corings was extended to comodules   
by El Kaoutit and G\'omez-Torrecillas in \cite{ElGo.com}, where
to any bimodule $_SP_A$ with $P_A$ finitely generated and projective, 
a coring $P^*\ot_S P$ was associated and it was shown that the canonical
map  
$$\tilde\mu_A:\Hom_A(P,A)\ot_SP\to \cC$$
is a coring morphism provided $P$ is also a right $\cC$-comodule and
$S=\End^\cC(P)$.
In \cite[18.25]{BW} such comodules $P$ were termed {\em Galois comodules}
provided $\tilde\mu_A$ was bijective, and it was proved 
in \cite[18.26]{BW} that this condition implies that 
the functors $\Hom_A(P,-)\ot_SP$ and $-\ot_A\cC$ from the right $A$-modules
to the right $\cC$-comodules are isomorphic.  

In a recent paper \cite{Brz:gal}, Brzezi\'nski further investigated
these Galois comodules and pointed out their relevance for 
{\em descent theory}, {\em vector bundles},  
and {\em non-commutative geometry}. Related questions are, for example, also 
considered by Caenepeel, De Groot and Vercruysse in \cite{CaGrVe}.
In this note we   concentrate on comodule properties and we want to free the
notion from the condition that $P_A$ has to be finitely generated and
projective.  This is done by taking the above mentioned isomorphism of
functors as definition. Although some symmetry is lost those properties which
to us seem to be essential, are preserved. From this point of view Galois
comodules are somehow similar to tilting (co)modules, or modules $M$
for which all $M$-generated modules are $M$-static: they all      
share the property that they are generators in their
respective categories provided they are flat over their endomorphism 
rings. Hence the presentation is partly motivated by the papers 
\cite{WiTil,WiStat} on tilting and static modules. Some results from 
\cite{Brz:gal} and \cite{CaGrVe} are obtained in a more general setting
(e.g., \ref{rem.Br}, \ref{rem.Ca}). 

Relative injectivity of comodules  
is of special interest in the context of our investigations and leads to
category equivalences. In particular, a strongly $(\cC,A)$-injective
(equivariantly injective)  Galois comodule $P$ that is 
finitely generated and projective as $A$-module, induces an equivalence
between  the category of comodules and the $\End^\cC(P)$-modules 
(see \ref{equival}). 

Given a commutative ring $R$, 
an {\em entwining structure} $(A,C,\psi)$ consists of an $R$-algebra $A$,
an $R$-coalgebra $C$, and an $R$-linear map $\psi: C\ot_RA\to A\ot_R C$    
satisfying certain compatibility conditions which ensure that 
$A\ot_R C$ allows for an $A$-coring structure (see \cite[Section 32]{BW}). 
If $A$ is a right $A\ot_RC$-comodule (equivalently, there exists 
a grouplike element in $A\ot_RC$), then $\Hom^{A\ot_RC}(A,A\ot_RC)\simeq A$ 
and $A$ is a Galois $A\ot_RC$-comodule if and only if $A$ is a $C$-Galois
extension over $\End^{A\ot_RC}(A)=A^{coC}$ (see \cite[34.10]{BW}).
 Hence a number of results on 
generalised Hopf Galois (coalgebra Galois) extensions in Schauenburg-Schneider
\cite{SchSch} can be seen as special cases of our results (see \ref{rem.Br}). 

The symmetry for Galois comodules that are finitely generated and projective
$A$-modules mentioned above can be maintained for comodules which are 
direct sums of comodules of this type. In this context the 
{\em infinite comatrix corings}, as introduced by El Kaoutit
and G\'omez Torrecillas in \cite{ElGo}, find a natural application
(Section 6).

\section{Preliminaries}

Throughout we will essentially follow the notation in \cite{BW}. For
convenience we  recall some basic notions.  

\begin{thm}{\bf Corings.} \em
Let $A$ be an associative ring with unit and $\cC$ an $A$-coring with 
coproduct and counit
$$\Delta:\cC\to \cC\ot_A\cC, \quad \vareps: \cC\to A.$$
Associated to this there are the right and  left  dual rings
$\cC^*=\Hom_A(\cC,A)$ and ${^*\cC}={_A\Hom(\cC,A)}$ with 
the convolution products.
\end{thm}

\pagebreak[3]

\begin{thm}{\bf Comodules.} \em
A right $A$-module $M$ is a right $\cC$-comodule provided there is an 
$A$-linear $\cC$-coaction 
 $$\varrho^M: M\to M\ot_A\cC, \mbox{ written as } 
 \varrho^M(m)= \sum m\sw{0}\ot m\sw{1} \mbox{ for } m\in M,$$ 
satisfying the coassociativity and counital condition. 

We denote the category of right $A$-modules by $\cM_A$ and the category 
of right $\cC$-comodules by $\cM^\cC$. The corresponding left versions are 
denoted by $_A\cM$ and ${^\cC\cM}$, respectively. 
The category $\cM^\cC$ is additive, has coproducts and  
  cokernels, and epimorphims are surjective maps. 
 The functor 
$-\ot_A\cC: \cM_A\to \cM^\cC$ is right adjoint to the forgetful functor
by the isomorphisms, for $M\in \cM^\cC$, $X\in \cM_A$,
        $$\varphi: \Hom^\cC  (M, X\ot_A \cC ) \to \Hom_A( M, X), \; 
         f \mapsto  (\id_{X}\ut  \ueps)\circ f, $$        
         with inverse map \;  
	$h\mapsto (h\ut  \id_{\cC})\circ \varrho^M$. 

Notice that for any monomorphism (injective map) $f:X\to Y$ in $\cM_A$, the
colinear map
$f\ot\id_\cC:X\ot_A\cC\to Y\ot_A\cC$ is a monomorphism in $\cM^\cC$
but need not be injective. In case $_A\cC$ is flat, monomorphisms in 
$\cM^\cC$ are injective maps and in this case $\cM^\cC$ is a Grothendieck
category (see \cite[18.14]{BW}).
\end{thm}

 \begin{thm}\Label{sigmaM}{\bf The subcategory $\sM$.} \em
Let $M\in \cM^\cC$. Homomorphic images of direct sums of copies of $M$
are called {\em $M$-generated} comodules. The full subcategory of
$\cM^\cC$, whose objects are subcomodules $K$ of $M$-generated
comodules $N$ (i.e. there is an injective colinear map $K\to N$), 
is denoted by $\sM$.
Notice that this does not imply that morphisms in $\sM$ have kernels 
unless $_A\cC$ is flat. 

Cokernels of morphisms $M^{(\Lambda')} \to M^{(\Lambda)}$, with any sets
$\Lambda', \Lambda$, are called {\em $M$-presented comodules}.
Notice that the image of the functor $-\ot_SM: \cM_S\to \cM^\cC$ lies  
in $\sM$, in fact, comodules of the form $X\ot_SM$ with $X\in \cM_S$
are $M$-presented. 
\end{thm} 

\begin{thm}\Label{cond.alpha}{\bf  The $\alpha$-condition.} \em 
 Defining the convolution product on $^*\cC={_A\Hom}(\cC,A)$ as in \cite{BW},
any right $\cC$-comodule $(M,\varrho^M)$ allows a left $^*\cC$-module structure
by putting  $ f\rhu m = (\id_M\ot f)\circ\varrho^M(m)$, for any $f\in {^*\cC}$,
$m\in M$. This yields a faithful functor  $\Phi: \cM^\cC \to {_{^*\cC}\cM}$ 
 which is a full embedding if and only if  the map
 $$\alpha_{K}: K\ot_A \cC  \to \Hom_A(^*\cC,K), \quad n\ot
c\mapsto [f\mapsto nf(c)],$$ 
is injective for any $K\in \cM_A$. This holds if and only if 
$_A\cC$ is locally projective and is called 
{\em left $\alpha$-condition} on $\cC$. 
In this case $\MC$ can be identified with
$\sigma[{_{^*\cC} \cC}]$, the full subcategory of ${_{^*\cC} \cM}$
whose objects are subgenerated by $\cC$. 

Symmetrically the {\em right $\alpha$-condition} is defined and if it holds
$\CM$ can be identified with the category $\sigma[\cC_{\cC^*}]$ where
${\cC^*}=\Hom_A(\cC,A)$.
\end{thm}

\begin{thm}\Label{mor.def}{\bf Morphism groups.} \em
 The comodule morphisms between $M,N\in \cM^\cC$ is characterised by the
exact sequence of $\Z$-modules 
 $$0\to \Hom^\cC
(M,N)\to\Hom_A(M,N)\stackrel{\gamma}{\lra}\Hom_A(M,N\ot_A \cC ),$$   where
$\gamma (f)=\varrho^N\circ f - (f\ut \id_{\cC})\circ \varrho^M$
\end{thm}

\begin{thm}\Label{cot.def}{\bf Cotensor product.} \em For two comodules 
$M\in \CM$ and $L\in \CM$ the cotensor product is define by the 
exact sequence of $\Z$-moduls 
$$\xymatrix{ 
 0\ar[r] & M \otC L \ar[r]  & 
    M \ot_A N\ar[r]^{\omega_{M ,L}\quad} & 
          M \ot_A   \cC\ot_A L }$$
where $\omega_{M,L}=\roM\ot\id_{L}-\id_{M} \ot {^L\varrho}$.
\end{thm}

\begin{thm}\Label{fp.fg}   
{\bf $M_A$ finitely generated projective.} \em
Let $M\in \cM^\cC$ such that 
$M_A$ is finitely generated and projective.
Then for $M^*=\Hom_A(M,A)$ the map
 $$ \varphi: \cC \ot_A M^* \to \Hom_A(M,\cC), \quad c\ot h\mapsto c\ot h(-)$$
is an isomorphism and induces a left $\cC$-comodule structure 
on $M^*$ (see \cite[19.19]{BW}).
With a dual basis $m_1,\ldots,m_n \in M$,
$\pi_1,\ldots,\pi_n \in M^*$, the 
the inverse map of $\varphi$ is given by sending 
$g\in M^*$ to $\sum_i g(m_i)\ot\pi_i$,   
and the coaction on $M^*$ is 
$$\begin{array}{c}
\varrho^{M^*}: M^* \to \cC\ot_A M^*, \quad
 g\mapsto(g\ot\id_\cC)\varrho^M\mapsto\sum_i(g\ot \id_\cC)\varrho^M(m_i)\ot
\pi_i.
\end{array}$$ 
  
There is a canonical anti-isomorphism between  $^\cC\End({M^*})$ and
$S=\End^\cC(M)$ and by this $M^*$ is a right $S$-module.  

For any $N\in \cM^\cC$, there exists an
isomorphism (natural in $M$) 
$$ N\otC M^*\stackrel{\simeq}{\lra} \Hom^\cC(M,N).$$ 

This follows from the proof of \cite[10.11]{BW}: 
With the defining sequences for $\Hom^\cC$ and $\otC$ we have the 
commutative diagram with exact rows 
$$\xymatrix{ 
 0\ar[r] & N \otC M^* \ar[r]\ar@{-->}[d] & 
    N \ot_A M^*\ar[r]^{\omega_{N ,M^*}\quad} \ar[d]^\simeq & 
          N \ot_A   \cC\ot_A  M^*\ar[d]^\simeq  \\ 
 0\ar[r] & \Hom^\cC(M,N )\ar[r] & \Hom_A(M,N )\ar[r]^{\gamma\quad\;} & 
\Hom_A(M,N \ot_A  \cC),   }$$ 
where $\omega_{N ,M^*}=\roN  \ut \id_{M^{*}} -  
\id_{N } \ut \varrho_{M^*}$  
and $\gamma (f) := \roN \circ f - (f\ut \id_{\cC})\circ \varrho^M$.
From this diagram lemmata imply the existence and bijectivity of the required
morphism.  
\end{thm}

Notice that this isomorphism is also proved in \cite[Proposition 4]{CaGrVe}.

\begin{thm}\Label{coint}{\bf Cointegrals.} \em
An $(A,A)$-bilinear map $\delta:\cC\ot_A\cC\to \cC$ is called a {\em cointegral
in $\cC$} if   
$$(\id_\cC \ot\delta)\circ (\Delta \ut \id_\cC  )= 
      (\delta\ot \id_\cC)\circ (\id_\cC  \ot\Delta).$$
Cointegrals are characterised by the fact that 
 for any $M\in \cM^\cC$, the map 
 $$\nu_M=(\id_M\ot\delta)\circ (\varrho^M\ot\id_\cC): M\ot_A \cC \to M$$ 
     is a comodule morphism, or by the 
corresponding property for left $\cC$-comodules. 

This follows from the proof of \cite[3.29]{BW}. In
\cite[Section 5]{CaGrVe} these maps are related to the counit for the adjoint
pair of functors $-\ot_A\cC$ and the forgetful functor.
For $R$-coalgebras $C$ over a commutative ring $R$ with $C_R$ locally
projective, a cointegral is precisely a $C^*$-balanced $R$-linear map $C\ot_R
C\to R$  (e.g., \cite[6.4]{BW}).
\end{thm}

\begin{thm}\Label{rel.inj}{\bf Relative injectivity.} \em 
 Let $M$ be a right $\cC$-comodule and $S= \End^\cC(M)$.  

$M$ is {\em $(\cC,A)$-injective} provided 
the structure map $\varrho^M:M\to M\ot_A\cC$ is split by a $\cC$-morphism 
$\lambda: M\ot_A\cC\to M$. 

We call $M$ {\em strongly $(\cC,A)$-injective} if this $\lambda$ is
$\cC$-co\-linear and $S$-linear.
Given a subring $B\subseteq S$,  $M$ is said to be
{\em $B$-strongly $(\cC,A)$-injective} if $\lambda$ is 
$\cC$-co\-linear and $B$-linear.
 
We call $M$ {\em fully $(\cC,A)$-injective} if there exists 
 a cointegral $\delta_M:\cC\ot_A\cC\to \cC$ such that 
$\roM$ is split by $(\id_M\ot \delta_M)\circ (\varrho^M\ot \id_\cC)$.

The notions for left $\cC$-comodules are defined symmetrically. 
\end{thm}

Obviously fully $(\cC,A)$-injective are strongly $(\cC,A)$-injective.
and for a 
$B$-strongly $(\cC,A)$-injective comodule $M$ 
and any $X\in \cM_B$, $X\ot_B M$ is $(\cC,A)$-injective.

For coalgebras $B$-strongly $(\cC,A)$-injective comodules are named 
{\em $B$-equivariantly $\cC$-injective}
  (see \cite[Definition 5.1]{SchSch}).
Cointegrals $\delta$ making  $M$ fully $(\cC,A)$-injective 
are said to be {\em $M$-normalized} in \cite[Proposition 5.1]{CaGrVe}.

The fact that under projectivity conditions comodule properties may be
considered as module properties has the following implication.

\begin{thm}\Label{C.strongly} {\bf $\cC$ strongly $(\cC,A)$-injective.}
 Assume $\cC_A$ to be locally projective. Then the following are equivalent:
\begin{blist}
\item $\cC$ is strongly $(\cC,A)$-injective;
\item $\cC$ is a coseparable coring.
\end{blist}
\end{thm}
\proof One implication is obvious. Recall that $\EndC(\cC)\simeq \cC^*$ 
and assume $\cC$ to be strongly $(\cC,A)$-injective
with a $\cC^*$-splitting right $\cC$-colinear map $\nu:\cC\ot_A\cC \to \cC$. 
By the right $\alpha$-condition this means that $\nu$ is also left
$\cC$-colinear and hence $\cC$ is coseparable.
\endproof

\begin{thm}\Label{full.prop}
{\bf Properties of fully $(\cC,A)$-injective comodules.}  Let 
$M\in \cM^\cC$ with $S=\End^\cC(M)$. 
\begin{zlist}
\item $M$ is fully $(\cC,A)$-injective if and only if 
   $$ (\id_M\ot \widetilde\delta_M) \circ\varrho^M =\id_M\;\mbox{  
  where }\;\widetilde\delta_M=\delta_M \circ \Delta: \cC\to A.$$ 
\item $\cC$ is a fully $(\cC,A)$-injective left (right) comodule 
      if and only if $\cC$ is a coseparable coring.
\item Let $M$ be fully $(\cC,A)$-injective. Then:
\begin{rlist}
 \item Every comodule in $\sM$ is fully $(\cC,A)$-injective.   
 \item If $M$ is a subgenerator in $\MC$ then $\cC$ is a coseparable coring. 
 \item For any subring $B\subset S$ and $X\in \cM_B$, $X\ot_B M$ is  
      fully $(\cC,A)$-injective.
 \item If $M_A$ is finitely generated and projective, then 
      $M^*$ is a fully $(\cC,A)$-injective left $\cC$-comodule. 
\end{rlist}
\end{zlist} 
\end{thm}
\proof
(1) The assertion follows from the equalities
  $$ (\id_M\ot \widetilde\delta_M) \circ\varrho^M =
(\id_M\ot\delta_M)\circ(\varrho^M\ot\id_\cC)\circ\varrho^M=\id_M.$$ 

(2) This is shown in \cite[26.1]{BW}.
   In this case $\delta_\cC = \vareps$.  
\smallskip

(3)(i) Obviously any direct sum $M^{(\Lambda)}$ is fully $(\cC,A)$-injective.
  For every $\cC$-comodule epimorphism $f:M\to N$ and $m\in M$,
   $$(\id_N\ot\delta_M)(f(m)\sw{0}\ot f(m)\sw{1})=
   (\id_N\ot\delta_M)(f(m\sw{0})\ot
m\sw{1})=f(m\sw{0}\delta_M(m\sw{1}))=f(m).$$ 
This proves that factor comodules 
of $M$ are fully $(\cC,A)$-injective. 
For subcomodules similar arguments apply.

(ii) This follows from (2) and (3)(i).
 
(iv) With the dual basis $(m_i,\pi_i)$ for $M$, the coaction of $\cC$ on
$g\in M^*$ is given by  
$ \varrho^{M^*}(g)=\sum_i (g\ot \id_\cC)\varrho^M(m_i)\ot\pi_i$, and 
$$\begin{array}{c}
\sum_i(\delta_M\ot \id_{M^*})(g({m_i}\sw{0}){m_i}\sw{1} \ot \pi_i) =
 \sum_i g({m_i}\sw{0}\delta_M({m_i}\sw{1}))\ot\pi_i=
  \sum_i g(m_i)\ot\pi_i = g,
\end{array}$$     
proving that $M^*$ is a fully $(\cC,A)$-injective comodule.
 This isomorphism is also proved in \cite[Proposition 5]{CaGrVe}.
\endproof

\begin{thm}\Label{split} {\bf Splitting of $\Hom$- and $\ot$-sequences.}
\begin{zlist}
\item Let $M\in {\cM^\cC}$, $B\subseteq \End^\cC(M)$ a subring, and
 assume $M$ to be $B$-strongly $(\cC,A)$-injective. Then 
 for any $N\in {\cM^\cC}$, the sequence  
$$\xymatrix{ 
 0\ar[r] & M \otC L \ar[r]  & 
    M \ot_A L\ar[r]^{\omega_{M,L}\quad} & 
    M \ot_A\cC\ot_A L}$$
is splitting in $_B\cM$, 
 and for any $L\in {^\cC\cM}$ the  sequence  
$$\xymatrix{ 
 0\ar[r] & \Hom^\cC(N,M )\ar[r] & \Hom_A(N,M )\ar[r]^{\gamma\quad\;} & 
   \Hom_A(N,M \ot_A  \cC)}$$ 
   is also splitting in $_B\cM$.
\item Let $L\in {^\cC\cM}$, $D\subseteq {^\cC\End}(L)$ a subring, 
   and assume $L$ to be $D$-strongly $(\cC,A)$-injective. Then
 for any $M\in {\cM^\cC}$, the sequence   
$$\xymatrix{ 
 0\ar[r] & M \otC L \ar[r]  & 
    M \ot_A L\ar[r]^{\omega_{M ,L}\quad} & 
          M \ot_A   \cC\ot_A  L }  $$
is splitting in $\cM_D$, and for any $K\in {^\cC\cM}$, the sequence
$$\xymatrix{ 
 0\ar[r] & \CHom(K,L)\ar[r] & \Hom_A(K,L)\ar[r]^{\gamma\quad\;} & 
\Hom_A(K,\cC\ot_AL) }$$ 
is also  splitting in $\cM_D$.
\end{zlist}  
\end{thm}
\proof  
(1) Let $\nu:M\ot_A\cC\to M$ be a comodule splitting of $\roM$. 

 As in the proof of \cite[21.5]{BW}(4) it is easy to see that  
$$\beta = (\nu\ot \id_L)\circ(\id_M\ot \roL):M\ot_A N\to M\otC N$$  
is an $\End^\cC(L)$-linear retraction. If $\nu$ is $B$-linear
then  obviously $\beta$ is a $B$-linear retraction and the proof
of \cite[21.5]{BW}(4) applies.

 For the second sequence 
we can follow the proof of \cite[3.18]{BW}.
The inclusion is split by
 $$\Hom_A(N,M)\to \Hom^\cC(N,M):f\mapsto \nu\circ(f\ot\id_\cC)\circ
\varrho^N,$$
 and $\gamma$ is split modulo $\Hom^\cC(N,M)$ by 
$$\Hom_A(N,M\ot_A\cC)\to \Hom_A(N,M), \quad g\mapsto\nu \circ g.$$
This is clearly a right $\EndC(N)$-linear splitting. 
If $\nu$ is $B$-linear, the splitting maps are also left
$B$-linear.  
\smallskip

(2) The assertions can be seen by symmetry.
\endproof

For convenience we list some

\begin{thm}\Label{ass.cot} {\bf Associativity conditions for the cotensor
  product.}
   Consider two comodules $M\in {\cM^\cC}$ and $L\in {^\cC\cM}$.
\begin{zlist}
\item For a subring $B\subset \End^\cC(M)$ and $X\in \cM_B$, 
   $$X\ot_B (M\ot^\cC L ) \simeq  (X\ot_B M)\ot^\cC L$$ 
 provided that 
  \begin{tabular}[t]{rl}
  {\rm (i)}& $X$ is a flat $B$-module, or \\
  {\rm (ii)} & $-\otC L$ is right exact
         or $L$ is $(\cC,A)$-injective, or\\
  {\rm (iii)}& $M$ is $B$-strongly $(\cC,A)$-injective.
\end{tabular}

\item For a subring $D\subset {^\cC\End}(L)$ and $Y\in {_D\cM}$,
   $$  M\ot^\cC (L \ot_D Y)\simeq (M\ot^\cC L )\ot_D Y $$ 
 provided that 
  \begin{tabular}[t]{rl}
   {\rm (i)}& $Y$ is a flat $D$-module, or \\
   {\rm (ii)} & $M\otC-$ is a right exact or $M$ is $(\cC,A)$-injective,
  or \\   
 {\rm (iii)}& $L$ is  $B$-strongly $(\cC,A)$-injective.
 \end{tabular}
\end{zlist}
\end{thm}
\proof
The conditions (i),(ii) are sufficient
  to imply the assertion by \cite[21.4 and 21.5]{BW}.
 The sufficiency of (iii)  follows from \ref{split}.
\endproof

\begin{thm}\Label{hom.cot}
{\bf Hom-tensor relation.}   For $M\in \MC$, $L\in \CM$ a subring 
$B\subseteq \EndC(M)$ and any right $B$-module $X$, there is a map
$$\psi_X:X\ot_B\Hom^\cC(N,M)\to \Hom^\cC(N,X\ot_B M),\quad h\ot x\mapsto
x\ot h(-),$$ 
and this is an isomorphism provided 
\begin{rlist} 
 \item $X$ is a flat $B$-module and $N_A$ is finitely presented, or
 \item $M$ is $B$-strongly $(\cC,A)$-injective and $N_A$ is finitely generated 
       and projective, or
 \item $N$ is projective in $\cM^\cC$ and $N_A$ is finitely generated.
\end{rlist}
\end{thm}
\proof
  Consider the commutative diagram with canonical maps 
  $$\xymatrix{ 
 0\ar[r] & X\ot_B \HomC (N,M) \ar[r] \ar@{-->}[d]   & 
     X\ot_B \Hom_A (N, M) \ar[r] \ar[d]  &
      X\ot_B  \Hom_A (N, M\ot_A\cC) \ar[d] \\
 0\ar[r] &   \HomC (N,X\ot_BM) \ar[r]& \Hom_A (N,X\ot_BM)\ar[r] &
      \Hom_A (N,X\ot_B M\ot_A\cC), } $$
where the bottom sequence is exact.
If (i) holds then the top sequence is also exact, and by \ref{split}, this is
also true if two holds. In both cases the two right vertical maps are
isomorphisms and hence the first one is also an isomorphism. 

Now assume (iii) and consider an exact sequence $F_1\to F_2\to X\to 0$ in
$_B\cM$ where $F_1$, $F_2$ are free $B$-modules.  
With $-\ot_B M$ and $\HomC(N,-)$ we construct the commutative diagram with 
exact rows 
 $$\xymatrix{ 
 F_1\ot_B\HomC(N,M)\ar[r]\ar[d]^{\psi_{F_1}} &
  F_2\ot_B\HomC(N,M)\ar[r]\ar[d]^{\psi_{F_2}}&
       X\ot_B\HomC(N,M)\ar[r]\ar[d]^{\psi_X} &0\\ 
 \HomC(N,F_1\ot_B M)\ar[r]  & \HomC(N,F_2\ot_B M)\ar[r] &
    \HomC(N,X\ot_B M)\ar[r] &0 }$$
where  $\psi_{F_1}$ and $\psi_{F_2}$ are isomorphisms (since $\HomC(N,-)$
commutes with direct sums) and hence ${\psi_X}$ is an isomorphism.

Notice that projectivity of 
the comodule $N$ implies projectivity of $N_A$ (see \cite[18.20]{BW}).
Hence $\HomC(N,M)\simeq M\otC N^*$ and $N^*$ is coflat. 
So the assertion also follows from \ref{ass.cot}(1)(ii).
\endproof

\section{Adjoint functors and static comodules}

\begin{thm}\Label{adjoint} {\bf Adjoint pair of functors.} \em
 For any right $\cC$-comodule $P$ with endomorphism ring $S=\End^\cC(P)$,
 the functors (see \cite[18.21]{BW})
$$ -\ot_S P:\cM_S\to  {\cM^\cC},\quad  {\Hom^\cC}(P,-):{\cM^\cC}\to\cM_S,$$
form an adjoint pair by the functorial isomorphism 
  (for $N\in {\cM^\cC}$ and $X\in \cM_S$),  
$$\Hom^\cC(X\ot_S P , N) \to \Hom_S(X,\Hom^\cC(P,N)),\;
   g\mapsto [x\mapsto g(x\ut-)]\,, $$ 
with inverse map\; $h \mapsto [x\ut p \mapsto h(x)(p)]$. 
Counit and unit of this adjunction are given by
$$\begin{array}{ll}    
  \mu_N: {\Hom^\cC}(P,N)\ot_S P \to N, \ &f\ut p \mapsto f(p),\\[+1mm]       
  \nu_X: X\to {\Hom^\cC}(P,X\ot_S P),  &  x \mapsto[p\mapsto x\ut p] ,  
\end{array} $$ 
and each of the following compositions of maps yield the identity,
$$  
\xymatrix{
\Hom^\cC(P,N )\ar[rr]^{\nu_{\Hom(P,N)}\qquad\quad}&&
     \Hom^\cC(P,\Hom^\cC(P,N)\ot_SP)           
\ar[rr]^{\quad\qquad\Hom(P,\mu_N)}&&\Hom^\cC(P,N), } $$
$$  
\xymatrix{
X\ot_S P \ar[rr]^{id\ot\nu_X\qquad\quad}& &
  \Hom^\cC(P,X\ot_S P)\ot_S P\ar[rr]^{\quad\qquad\mu_{X\ot P}}&& X\ot_SP.} $$  

A $\cC$-comodule $N$ is called {\em $P$-static} if $\mu_N$ is an isomorphism,
% and the class of all $P$-static $\cC$-modules is denoted by $\St(P)$. 
and an $S$-module $X$ is called {\em $P$-adstatic} if $\nu_X$ is an
isomorphism. 
% and we denote by $\Ad(P)$  the class of all $P$-adstatic $S$-modules. 
Cleary $P$ is $P$-static and this is also true for 
direct sums of copies of $P$ since, for any index set $\Lambda$,
$$\mu_{P^{(\Lambda)}}: \Hom^\cC(P,P^{(\Lambda)})\ot_SP\to P^{(\Lambda)},$$
is a comodule isomorphism with inverse map $(p_\lambda)_\Lambda \mapsto
\sum_\Lambda \epsilon_\lambda \ot p_\lambda$, 
where $\epsilon_\lambda$ and
$\pi_\lambda$ denote the canonical inclusions and projections of the coproduct
$P^{(\Lambda)}$.
\end{thm}

\begin{thm}\Label{pro.prop}{\bf $P_A$ finitely generated and projective.} 
Let $P\in \cM^\cC$ with 
$P_A$ finitely generated and projective and $S=\End^\cC(P)$.  
\begin{zlist}
\item For any $N\in \cM^\cC$, $\Hom^\cC(P,N)\simeq N\ot^\cC P^*$,
 in particular  $S\simeq P\ot^\cC P^*$.

\item A module $X\in \cM_S$ is $P$-adstatic, provided 
       $$(X\ot_{S}P)\ot^\cC P^*\simeq X\ot_{S}(P\ot^\cC P^*).$$
\item 
\begin{rlist}
\item Every flat $X\in \cM_S$ is $P$-adstatic.
\item If $P^*$ is coflat or $(\cC,A)$-injective, or $P$ is 
      strongly $(\cC,A)$-injective, 
      then every $X\in \cM_S$ is $P$-adstatic.
\end{rlist}  
\item For a subring $B\subset S$ and  $Y\in \cM_{B}$,
      $$\Hom^\cC(P,Y \ot_{B} P)\simeq Y\ot_{B} S$$
      provided $Y_B$ is flat, or  
      $P$ is $B$-strongly $(\cC,A)$-injective, or $P^*$
     is coflat or $(\cC,A)$-injective. 
\end{zlist}
\end{thm}
\proof 
(1) This is shown by the proof of \cite[21.8]{BW}.
\smallskip

(2) Recall that $X\in \cM_S$ is $P$-adstatic if 
    $\nu_X: X \to  \Hom^\cC(P, X\ot_S P)$
is an isomorphism. 
Under the given condition, (1) implies 
  $$\Hom^\cC(P, X\ot_S P)\simeq X\ot_S (P\ot^\cC P^*)\simeq X .$$

(3) As shown in \ref{ass.cot}, each of the conditions implies the
isomorphism required. 
\smallskip

(4) From (1) we get $\Hom^\cC(P,Y\ot_B P)\simeq (Y\ot_{B}P)\ot^\cC P^*$,
and by \ref{ass.cot}, under each of the conditions required,
$$(Y\ot_{B}P)\ot^\cC P^*\simeq Y\ot_{B}(P\ot^\cC P^*)\simeq Y\ot_{B}S.$$
\endproof

\begin{thm}{\bf $P$ as generator  in $\cM^\cC$.} \em
Recall that $P$ is a generator in $\cM^\cC$ if and only if 
 the functor $\Hom^\cC(P,-):\cM^\cC\to \cM_S$ is faithful and that faithful
functors reflect epimorphisms (e.g. \cite[11.3]{WiBu}). 
Since, for any $N\in \cM^\cC$,
$$  \Hom^\cC(P,\mu_N): \Hom^\cC(P,\Hom^\cC(P,N)\ot_SP)\to \Hom^\cC(P,N)$$ 
is an epimorphism (surjective), we conclude that  
$\mu_N$ is an epimorphism (that is surjective) in $\cM^\cC$ provided
$P$ is a generator in $\cM^\cC$. 
Taking $\Lambda = \Hom^\cC(P,N)$, the canonical epimorphism 
 $\varphi_N:P^{(\Lambda)}\to N$ remains an epimorphism under $\Hom^\cC(P,-)$.

Now assume $_A\cC$ to be flat. Then $K=\Ke\, \varphi$ is a comodule, and 
we have the commutative diagram with exact rows
$$\xymatrix{
 & \Hom^\cC(P,K)\ot_SP\ar[r]\ar[d]^{\mu_K} & \Hom^\cC(P,P^{(\Lambda)})\ot_SP
   \ar[d]^\simeq  \ar[r] & \Hom^\cC(P,N)\ot_SP \ar[d]^{\mu_N}\ar[r] & 0 \\
  0\ar[r] & K\ar[r] & P^{(\Lambda)} \ar[r] & N \ar[r] & 0,}$$
where ${\mu_K}$ is surjective. By diagram lemmata this implies that 
$\mu_N$ is injective (hence an isomorphism). 
\end{thm}

\pagebreak[3] 

\begin{thm}\Label{gen}{\bf Properties of generators.}
Assume $_A\cC$ to be flat.
\begin{zlist}
\item If $P$ generates the (finitely generated) subcomodules of $P^{(\N)}$,
     then, for every $P$-generated $\cC$-comodule $L$, $\mu_L$ is an
     isomorphism, and for every $N\in \cM^\cC$, $\mu_N$ is injective.
\item $P$ is a generator in $\cM^\cC$ if and only if $\mu_N$ is an
      isomorphism for any $N\in\cM^\cC$, i.e., every right $\cC$-comodule is
$P$-static. \end{zlist} 
 \end{thm}

\proof (1) Clearly the condition implies that $P$ generates the subcomodules  
of any direct sum of copies of $P$ and  bijectivity of $\mu_L$
 follows from the considerations above. 
The image of $\mu_N$ is the trace $\Tr(P,N)$ of $P$ in $N$ (sum of all 
$P$-generated subcomodules) and $\Hom^\cC(P,N)= \Hom^\cC(P,\Tr(P,N))$.
Since $\mu_{\Tr(P,N)}$ is bijective $\mu_N$ has to be injective.

(2) is a special case of (1). 
\endproof

Semisimple right comodules $P$ are defined by the fact that any monomorphism
$U\to P$ is a coretraction, that is subcomdules are direct summands. 
If $_A\cC$ is flat this is equivalent to  
$P$ being a (direct) sum of simple subcomodules and then any direct sum 
of copies of $P$ is semisimple.  

\begin{thm}\Label{semis}{\bf Semisimple comodules.}
Let $_A\cC$ be flat, $P\in \cM^\cC$ and $S=\End^\cC(P)$.
\begin{zlist}
\item The following are equivalent: 
 \begin{blist}
 \item $P$ is semisimple;
 \item for any set $\Lambda$, $\End^\cC(P^{(\Lambda)})$ is a von Neumann
       regular ring and, for any $N\in\cM^\cC$,      
      $\mu_N:\Hom^\cC(P,N)\ot_SP\to N$ is injective;
 \item for any set $\Lambda$, $\End^\cC(P^{(\Lambda)})$ is a regular ring, and 
      for any $L\in \sP$, \\
 $\mu_L:\Hom^\cC(P,L)\ot_SP\to L$ is an isomorphism.
 \end{blist}

\item If $P$ is finitely generated (in $\cM^\cC$), then the following are
         equivalent:
\begin{blist}
\item $P$ is semisimple;
\item $S$ is a right (left) semisimple ring, and  for any $N\in\cM^\cC$, \\  
         $\mu_N:\Hom^\cC(P,N)\ot_SP\to N$ is injective;
\item $S$ is a right (left) semisimple ring, and 
      for any $L\in \sP$, \\
 $\mu_L:\Hom^\cC(P,L)\ot_SP\to L$ is an isomorphism. 
\end{blist}

\item The following are equivalent:
\begin{blist}
\item $P$ is simple;
\item $S$ is a division ring, and  for any $N\in\cM^\cC$, 
      $\mu_N:\Hom^\cC(P,N)\ot_SP\to N$ is injective;
\item $S$ is a division ring, and 
      for any $L\in \sP$, $\mu_L:\Hom^\cC(P,L)\ot_SP\to L$ is an isomorphism.
\end{blist}
\end{zlist}
\end{thm}
\proof (1)
(a)$\Ra$(b)$\LRa$(c) For any $s\in S$, the image and the kernel are  
 direct summands in $P$. This implies that $S$ is von Neumann regular  
(e.g., \cite[37.7]{WiBu}). Since $P^{(\Lambda)}$ is also semisimple
the same argument shows that $\End^\cC(P^{(\Lambda)})$ is von Neumann regular.
Since $P$ generates all submodules of any $P^{(\Lambda)}$ the remaining
assertions follow from \ref{gen}(1).
\smallskip

(b)$\Ra$(a) Let $N\subset P$ be any subcomodule and construct the 
commutative diagram
$$\xymatrix{
   0\ar[r] & \Hom^\cC(P,N)\ot_S P \ar[r]\ar[d]^{\mu_N}& 
   \Hom^\cC(P,P)\ot_S P \ar[d]^\simeq \ar[r] &
    \Hom^\cC(P,P/N)\ot_S P \ar[d]^{\mu_{P/N}}  \\
  0\ar[r] & N \ar[r] & P \ar[r] &  P/N  \ar[r] & 0 ,} $$
in which the top row is exact by regularity of $S$ ($_SP$ is flat).
Clearly $\mu_{P/N}$ is an epimorphism and is injective by assumption.
This implies that ${\mu_N}$ is an epimorphism and hence $N$ is $P$-generated.
So there is some epimorphism $h:P^{(\Lambda)}\to N$. Considering 
$h$ as an endomorphism of $P^{(\Lambda)}$, the fact that  
$\End^\cC(P^{(\Lambda)})$ is regular implies that the image of $h$ is a direct 
summand in $P^{(\Lambda)}$ and hence in $P$ (see \cite[37.7]{WiBu}). 
This shows that $P$ is semisimple.
\smallskip

(2) The endomorphism ring of a finite direct sum of simple comodules is 
right (left) semisimple and this implies that  $\End^\cC(P^{(\Lambda)})$ is
von Neumann regular. Hence the proof of (1) applies. 
\smallskip

(3) By Schur's Lemma the endomorphism ring of a simple comodule is 
a division ring and again the proof of (1) applies. 
\endproof

Note that assertion (3) is also proved in \cite[Theorem 3.1]{Brz:gal} 
\smallskip

If $_A\cC$ is flat, a generator in $\cM^\cC$ is characterized by the fact 
that all comodules are $P$-static. 
This suggests the study of comodules $P$ by the classes of $P$-static 
modules. Transferring observations from module theory
we may consider the following cases: 

\begin{thm}\Label{stat} {\bf Some classes $P$-static.} \em 
 Consider the following conditions for $P\in \cM^\cC$:
\begin{zlist}
\item All comodules in $\cM^\cC$ are $P$-static. 
\item The class of $P$-generated modules is $P$-static.
\item The class of $P$-presented comodules is $P$-static.
\item The class of injective comodules in $\cM^\cC$ is $P$-static.
\item The class of $(\cC,A)$-injective comodules in $\cM^\cC$ is $P$-static.
\end{zlist}
\end{thm}

The first case was handled in \ref{gen} for $_A\cC$ flat.  In module
categories the second case  describes an important property of {\em
self-tilting} modules; for those an additional projectivity condition is
required (see  \cite[4.2]{WiTil},  \cite[4.4]{WiStat}). The
third case generalises  tilting modules (see \cite[4.3]{WiStat}).
For a module $P$, the corresponding property (4) essentially means that 
all $P$-injective modules in $\sigma[P]$ are $P$-static and 
- if $P$ is a balanced bimodule - this can be seen as descending chain
condition on certain matrix subgroups of $P$ 
(see \cite[5.4]{WiTil}, \cite{Zim}). 
In all these cases the functor $\Hom^\cC(P,-)$ induces equivalences between 
the $P$-static classes and the corresponding adstatic classes.
Properties of these classes correspond to properties of the module $P$.
For example, if the class of $P$-adstatic comodules
is closed under infinite coproducts, then 
$P$ has to be self-small, i.e., $\Hom^\cC(P,P^{(\Lambda)})\simeq
\Hom^\cC(P,P)^{(\Lambda)}$. 

If $_A\cC$ is flat, monomorphisms 
are injective maps and kernels exist in $\cM^\cC$, and hence most of the
proofs  for module categories can be transferred to $\cM^\cC$. In particular,
if $_A\cC$ is locally projective, $\cM^\cC$ can be identified with 
$\sigma[{_{^*\cC}\cC}]$ and the results mentioned immediately apply to 
comodules. Without such restrictions all the properties listed are also of
interest and deserve to be investigated elsewhere.
Here we will investigate the comodules characterised by the condition
required in (5).

\section{Galois comodules}

Throughout this section let $\cC$ be an $A$-coring,
$P\in \cM^\cC$ and put $S=\End^\cC(P)$, $T=\End_A(P)$. 
We consider the relationship between the two functors
 $$-\ot_A\cC \mbox{ and } \Hom_A(P,-)\ot_S P :\cM_A\to \cM^\cC.$$    
 
\begin{thm}\Label{gal.comod}{\bf Galois comodules.}  
{\em We call $P$ a} Galois comodule {\em if the following
equivalent conditions hold:}
  \begin{blist}
\item  The functors $-\ot_A\cC$ and $\Hom_A(P,-)\ot_S P$ are isomorphic;
\item  $\Hom_A(P,-)\ot_SP$ is right adjoint to the forgetful 
    functor $\cM^\cC\to \cM_A$, that is, for $K\in \cM_A$ and $M\in \cM^\cC$,
    there is a (bifunctorial) isomorphism
    $$\Hom^\cC(M, \Hom_A(P,K)\ot_S P) \to \Hom_A(M,K);$$
\item for any $K\in \cM_A$ there is a functorial isomorphism of comodules
    $$\tilde\mu_K: \Hom_A(P,K)\ot_S P\to K\ot_A \cC, \;g\ot p\mapsto (g\ot
     \id_\cC)\varrho^P(p);$$
\item every $(\cC,A)$-injective $N\in \cM^\cC$ is $P$-static, i.e.,
    $$\mu_N:\Hom^\cC(P,N)\ot_S P\to N, \; f\ot p\mapsto f(p),$$
   is an isomorphism (in $\cM^\cC$).
\end{blist}
\end{thm}
\proof We prove the equivalence of the conditions.   

 (a)$\LRa$(b) is clear since  $-\ot_A\cC$ is
right adjoint to the given forgetful functor and right adjoints are unique up
to functorial isomorphisms.
\smallskip

(b)$\LRa$(c) Both functors $\Hom_A(P,-)\ot_SP$ and $-\ot_A\cC$ are 
adjoints of the same forgetful functor and hence they are isomorphic. 
\smallskip

(c)$\Ra$(d) (see proof of \cite[18.26]{BW})
 Assume  $N\in \cM^\cC$ to be $(\cC,A)$-injective. Then, by
\cite[18.18]{BW}, the canonical sequence  
$$\xymatrix{ 
 0 \ar[r]&  \Hom^\cC(P,N)\ar[r]^i& \Hom_A(P,N)\ar[r]^{\gamma\qquad} & 
           \Hom_A(P,N\ot_A\cC)}$$ 
is (split and hence) pure in $\cM_S$, where 
$\gamma(f)=\varrho^N\circ f-(f\ut\id_\cC)\circ \varrho^P$. 
Hence tensoring with $_SP$ yields the commutative  diagram with exact rows, 
$$\xymatrix{  
  0\ar[r] & \Hom^\cC(P ,  N) \ot_S  P \ar[r] \ar[d]^{\mu_N}  & 
       \Hom_A(P,N) \ot_S  P \ar[r] \ar[d]^{\tilde\mu_N}  &   
   \Hom_A(P , N\ot_A\cC)\ot_S P \ar[d]^{\tilde\mu_{N\ot C}} \\ 
  0\ar[r] & N \ar[r] & N\ot_A\cC \ar[r] & N\ot_A\cC\ot_A\cC \, , 
} $$ 
where  the $\tilde\mu$'s are isomorphisms and so is $\mu_N$. 
\smallskip

(d)$\Ra$(c)  
 Since  $K\ot_A\cC$ is $(\cC,A)$-injective the assertion follows 
 from the commutative diagram of right $\cC$-comodule maps 
$$\xymatrix{ 
  \Hom^\cC(P,K\ot_A\cC)\ot_S P\ar[r]^{\qquad\qquad\mu_{K\ot C}}\ar[d]_\simeq & 
K\ot_A \cC \ar[d]^=\\      
\Hom_A(P,K)\ot_S P\ar[r]^{\qquad\tilde\mu_K} & K\ot_A\cC } 
 \;  \xymatrix{ 
   f\ot p \ar@{|->}[r] \ar@{|->}[d]& f(p)\ar@{|->}[d]^= \\      
 (\id\ot \eC)\!\circ\! f \ot p \ar@{|->}[r] &  
  \sum (\id\ot \eC)\circ f(p\sw{0})\ot p\sw{1}. } $$
\endproof

\begin{thm}\Label{prop.gal}{\bf Properties of Galois comodules.}
 Let $P\in \cM^\cC$ be a Galois comodule. Then:
\begin{zlist}
\item For any $(\cC,A)$-injective $N\in \cM^\cC$, there is an isomorphism
  $$\nu_{\Hom^\cC(P,N)}:\Hom^\cC(P,N)\to\Hom^\cC(P,\Hom^\cC(P,N)\ot_SP),$$  
  that is,  $\Hom^\cC(P,N)$ is $P$-adstatic.
\item For any $K\in \cM_A$, there is an isomorphism
$$\nu_{\Hom_A(P,K)}:\Hom_A(P,K)\to \Hom^\cC(P,\Hom_A(P,K)\ot_SP).$$  

\item There are right $\cC$-comodule isomorphisms 
   $$\Hom^\cC(P,\cC)\ot_SP\simeq \cC \simeq  \Hom_A(P,A)\ot_SP.$$
   
\item Since $T=\End^\cC(P, P\ot_A\cC)$, there is a $T$-linear isomorphism
   $$\begin{array}{l}
T\ot_S P \to P\ot_A\cC,\quad t\ot p \mapsto (t\ot \id_\cC)\varrho^P(p),
\; \mbox{ and } \\[+1mm] 
 P^*\ot_T P \ot_A\cC \simeq P^*\ot_T T\ot_S P \simeq P^*\ot_S P \simeq\cC.
\end{array}$$

\item For any $K\in \cM_A$ and index set $\Lambda$,
   $$\Hom^\cC(P, (K\ot_A\cC)^{\Lambda})\ot_SP\simeq
     \Hom_A(P,K)^{\Lambda}\ot_SP\simeq 
        K^{\Lambda}\ot_A\cC.$$
\item There are isomorphisms 
 $$ \begin{array}{l}
  \Hom_A(\cC,A)\simeq \Hom_A(P^*\ot_SP,A)\simeq \End_S(P^*),\\[+2mm]  
  \Hom^\cC(\cC,P)\simeq \Hom^\cC(P^*\ot_SP,P)\simeq \Hom_S(P^*,S),
\mbox{ and }\\[+2mm] 
 \Hom^\cC(P\ot_A\cC,P)\simeq \Hom^\cC(T\ot_SP,P)\simeq \Hom_S(T,S).
\end{array}$$ 
\end{zlist}
\end{thm}
\proof (1), (2) follow from the fact that the composition  
  $\Hom^\cC(P,\mu_N)\circ \nu_{\Hom^\cC(P,N)}$ yields the identity.

(3),(4) Put $N=\cC$ or $N=P\ot_A\cC$ in the characterising relations.  

(5) This follows from the fact that the product of $\Lambda$ copies
of $K\ot_A\cC$ in $\cM^\cC$ is iso\-morphic to $K^{\Lambda}\ot_A\cC$.

(6) Apply isomorphisms from (3),(4) and properties of adjoint functors 
(see \ref{adjoint}).
\endproof

\begin{thm}\Label{CA-inj} {\bf $(\cC,A)$-injective modules.}%{full.gal}
 Let $P$ be a Galois comodule.
\begin{zlist}
\item For $N\in \cM^\cC$ the following are equivalent:
\begin{blist}
\item $N$ is $(\cC,A)$-injective;
\item $\Hom^\cC(P,\varrho^N):\Hom^\cC(P,N)\to \Hom^\cC(P,N\ot_A\cC)$ 
      is a coretraction in $\cM_S$. 
\end{blist}
\item For $P$ the following are equivalent:
\begin{blist}
\item $P$ is $(\cC,A)$-injective;
\item the inclusion $i:S\to T$ is split by a right $S$-linear map.
\end{blist}
\item For $P$ the following are equivalent:
  \begin{blist}
  \item $P$ is  strongly $(\cC,A)$-injective;
  \item the inclusion $i:S\to T$ is split by a $(S,S)$-bilinear map.
  \end{blist}
In this case every $P$-static comodule is  $(\cC,A)$-injective.
\item For $P$ the following are equivalent:
\begin{blist} 
\item $P$ is fully $(\cC,A)$-injective; 
\item $\cC$ is a coseparable $A$-coring.
\end{blist}
In this case every comodule in $\cM$ is $P$-static and
   fully $(\cC,A)$-injective.
\end{zlist}
\end{thm}
\proof (1) By (a), $\varrho^N$ splits in $\cM^\cC$ and hence 
$\Hom^\cC(P,\varrho^N)$ splits in $\cM_S$. In turn, (b) yields a splitting 
of $\varrho^N$ by tensoring with $-\ot_SP$.
\smallskip

(2) follow from (1) since $T\simeq \Hom^\cC(P,P\ot_A\cC)$ as right $S$-module.
\smallskip 
 
(3) By \ref{split}, the inclusion $\HomC(P,P)\to \HomC(P,P\ot_A\cC)\simeq T$
is split as $(S,S)$-bimodule.

If $P$ is strongly $(\cC,A)$-injective, then for any $X\in \cM_S$, the tensor 
product is $(\cC,A)$-injective. So in particular $P$-static comodules
are $(\cC,A)$-injective. 
\smallskip 

(4)(a)$\Ra$(b) Since $P$ is a subgenerator this follows from \ref{full.prop}.

 (b)$\Ra$(a) Over a coseparable coring all comodules are
fully$(\cC,A)$-injective. 
\endproof 
 
Notice that the assertions (2) and (3) in \ref{CA-inj} are shown in 
\cite[Theorem 7.2]{Brz:gal} for f.g. projective $A$-modules.  
The arguments in \cite{Brz:gal} can also be adapted to general 
Galois comodules.

\begin{thm}\Label{gal.alpha}{\bf Galois comodules under the
$\alpha$-condition.} \em 
If $\cC$ satifies the left $\alpha$-condition, $\MC$ can be identified 
with the $^*\cC$-module category $\sigma[{_{^*\cC} \cC}]$ 
(see \ref{cond.alpha})
and Galois comodules may be explained in these terms.

By the ring anti-morphism $A\to {^*\cC}$
(see \cite[17.7]{BW}) any left ${^*\cC}$-module has a right
$A$-module structure. It follows from the functorial isomorphisms on
$_{^*\cC}\cM$ for $K\in \cM_A$,
  $$\Hom_A(-,K)\simeq \Hom_A({^*\cC}\ot_{^*\cC}-,K)
     \simeq {_{^*\cC}\Hom}(-, \Hom_A({^*\cC},K)),$$
that $\Hom_A({^*\cC},K)$ is $({^*\cC},A)$-injective, that is,
injective with respect to short exact sequences in $_{^*\cC}\cM$ 
which split in $\cM_A$. 
Moreover, since the canonical map
$$\gamma_K: K\to \Hom_A({^*\cC},K),\; k\mapsto [f\mapsto f k],$$   
is $A$-split by $f\mapsto f(\varepsilon)$, it follows that a left 
${^*\cC}$-module $K$ is $({^*\cC},A)$-injective if and only  if $\gamma_K$
splits in $_{^*\cC}\cM$.  
For any $P\in \cM^\cC$ and $K\in \cM_A$, there are morphisms 
$$\begin{array}{rcl}
 \Hom_A(P,K)\simeq \Hom^\cC(P,K\ot_A\cC)&\stackrel{i}{\lra}
 &{_{^*\cC}\Hom} (P,K\ot_A\cC)\\    
&  \stackrel{\Hom(P,\alpha)}{\lra}&{_{^*\cC}\Hom}(P,\Hom_A({^*\cC},K))
   \simeq \Hom_A(P,K),
\end{array}$$
where $i$ is the inclusion and $\alpha$ is the canonical map from
\ref{cond.alpha}.  It is straightforward to prove that the
composition of these maps yields the identity on  $\Hom_A(P,K)$.
Hence injectivity of $\alpha$ implies that $\Hom(P,\alpha)$ is an isomorphism 
and leads to the following statement.
\end{thm}

\begin{thm}\Label{prop.alpha}{\bf Proposition.}
Let $P\in \cM^\cC$ be a Galois comodule, assume $\cC$ to satisfy 
the $\alpha$-condition and put $S=\End^\cC(P)={_{^*\cC}\End}(P)$. Then
for any $K\in \cM_A$,
$${_{^*\cC}\Hom}(P,\Hom_A({^*\cC},K))\ot_SP
  \simeq  \Hom^\cC (P,K\ot_A\cC)\ot_SP \simeq  K\ot_A \cC,$$
 implying  $K\ot_A \cC \simeq \Tr(P,\Hom_A({^*\cC},K))$.
\end{thm}

\proof Combine the observations above with isomorphisms for Galois comodules.
Notice that  $\Hom_A({^*\cC},K)$ need not be a $\cC$-comodule but the
trace of $P$ yields a ${^*\cC}$-submodule lying in $\cM^\cC$.
The last isomorphism is a special case of
the  corresponding observation for modules in \cite[20.4]{WiBim}).
\endproof

\begin{thm}{\bf Semisimple base ring.} \em
If the ring $A$ is left semisimple (artinian semisimple), then 
all $A$-modules are projective and injective and $(\cC,A)$-injective comodules
are in fact $\cC$-injective. Moreover, the $\alpha$-condition is satisfied
and $\cM^\cC$ corresponds to the category $\sigma[{_{^*\cC}\cC}]$. 
In this case Galois comodules are just the comodules $P$ for which all 
injectives in $\sigma[{_{^*\cC}\cC}]$ are $P$-static. Such modules were
considered in \cite{WiStat}. 
\end{thm}

\begin{thm}{\bf Remarks.} \em
The ideas outlined in \ref{gal.alpha} can be used as guideline  to study
modules $M$ of Galois type for module categories over ring extensions
$B\to A$ by the condition that all $(A,B)$-injective $A$-modules
are $M$-static.  
\end{thm}

Notice that so far we did not make any assumptions neither on the $A$-module
nor on the $S$-module structure of $P$. Of course properties of this type 
influence the behaviour of Galois comodules and we look at the $S$-module
structure first. 

\begin{thm}\Label{S.prop} {\bf Module properties of $_SP$.} 
 Let $P\in \cM^\cC$ be a Galois comodule.  
\begin{zlist}
\item If $_SP$ is finitely generated, then $_A\cC$ is finitely generated.
\item If $_SP$ is Mittag-Leffler, then $_A\cC$ is Mittag-Leffler.
\item If $_SP$ is finitely presented, then $_A\cC$ is finitely presented.
\item If $_SP$ is projective, then $_A\cC$ is projective. 
\item If $_TP$ is finitely generated and
      $_SP$ is locally projective, then $_A\cC$ is locally projective. 
\item  If $_SP$ is flat, then $_A\cC$ is flat and $P$ is a generator in
      $\cM^\cC$.
\item If $_SP$ is faithfully flat, then $_A\cC$ is flat and $P$ is a
projective       generator in $\cM^\cC$.
\end{zlist}
\end{thm}
\proof
(1),(2),(3) Putting $K=A$ in \ref{prop.gal}(6) we have the commutative diagram
 $$\xymatrix{
   \Hom_A(P,A)^{\Lambda} \ot_SP\ar[r]^{\qquad\simeq} \ar[d]^{\varphi_P} &
  A^{\Lambda}\ot_A\cC \ar[d]^{\varphi_\cC} \\
   (\Hom_A(P,A)\ot_SP)^{\Lambda}\ar[r]^{\qquad\quad\simeq} & \cC^{\Lambda}, }$$
where the $\varphi$'s denote the canonical maps. Then (e.g. \cite[12.9]{WiBu})
\smallskip

\begin{tabular}{rcccccl} 
$_SP$ is fin. gen.& $\Ra$ & ${\varphi_P}$ surjective &
   $\Ra$ & ${\varphi_\cC}$ surjective & $\LRa$ & $_A\cC$ fin. gen., \\
$_SP$ is ML & $\Ra$ & ${\varphi_P}$ injective &
   $\Ra$ & ${\varphi_\cC}$ injective & $\LRa$ & $_A\cC$ ML, \\
$_SP$ is fin. pres.& $\Ra$ & ${\varphi_P}$ bijective &
   $\Ra$ & ${\varphi_\cC}$ bijective & $\LRa$ & $_A\cC$ fin. pres.. \\
\end{tabular}
\\[+1mm]
 Recall that by definition $\cC$ is Mittag-Leffler (ML) if ${\varphi_\cC}$ 
is injective.
\smallskip

(4) Let $_SP$ be  projective. Then $T\ot_SP\simeq P\ot_A\cC$
is projective as left $T$-module. Consider any epimorphism $F\to \cC$ 
where $F$ is a free module in $_A\cM$. 
Then $\id_P\ot f$ is a splitting epimorphism in $_T\cM$, 
and in the commutative diagram with exact rows 
$$\xymatrix{
  P^*\ot_T P\ot_A F \ar[r]^{\id\ot\id\ot f}\ar[d]&
    P^*\ot_T P\ot_A \cC \ar[r]\ar[d]^\simeq &0 \\
   F  \ar[r]^f & \cC \ar[r] & 0 , }$$
where the first vertical map is the evaluation and the right 
 isomorphism is from \ref{prop.gal}(4), 
the top row is splitting in $_A\cM$ and hence 
$f$ also splits showing that $_A\cC$ is projective. 
\smallskip

(5) Let $_SP$ be locally projective.  
To check local projectivity of $_A\cC$ consider the diagram in $_A\cM$
 with $k\in \N$ and exact bottom row, 
 $$ \xymatrix{
    &\; A^k\;\ar[r]^{i} & \; \cC \ar[d]^{g}   & \\ 
   &\; L\ar[r]^{f} &  \; N\ar[r] \, & 0. } $$ 
Applying $P\ot_A-$ we obtain the diagram 
$$ \xymatrix{
    &\; P^k\;\ar[r]^{\id_P\ot i\quad} & \; P\ot_A\cC \ar[d]^{\id\ot g}   & \\ 
   &\; P\ot_A L\ar[r]^{\id\ot f} &  \; P\ot_A N\ar[r] \, & 0. } $$ 
Since  $P\ot_A\cC\simeq T\ot_S P$ is a locally projective $T$-module 
(by \cite[42.11]{BW})  and  
$_TP^k$ is finitely generated by assumption, there is some
$T$-morphism $h:P\ot_A\cC\to P\ot_A L$ with 
$$ ({\id\ot f})\circ h\circ ({\id_P\ot i}) = \id\ot g.$$ 
Applying $P^*\ot_T-$ and the evaluation map we obtain 
  $ f \circ (\id_{P^*}\ot h)\circ i = g$. This shows that 
$_A\cC$ is locally projective. 
\smallskip

(6) We have $-\ot_A\cC\simeq \Hom_A(P,-)\ot_S P$. Clearly $\Hom_A(P,-)$
(always) preserves injective maps. If $_SP$ is flat then $-\ot_S P$
also preserves injectivity of morphisms and hence $-\ot_A\cC$   
preserves injective maps, i.e., $_A\cC$ is flat.  

For any $M\in \cM^\cC$ we have an exact sequence of comodules
 $$\xymatrix{
0\ar[r]& M\ar[r]^{\varrho^M\quad}& M\ot_A\cC \ar[r] & M\ot_A\cC\ot_A\cC}.$$ 
 By left exactness of $\Hom^\cC(P,-)$
and $-\ot_S P$, we obtain the exact commutative diagram
$$\xymatrix{ 
 0 \ar[r]&  \Hom^\cC(P,M)\ot_S P\ar[r] \ar[d] & 
\Hom_A(P,M)\ot_S P \ar[d]^\simeq \ar[r] &\Hom_A(P,M\ot_A\cC)\ot_S P
\ar[d]^\simeq     \\ 
  0 \ar[r]& M \ar[r]^{\varrho^M}& M\ot_A\cC\ar[r] & M\ot_A\cC \ot_A\cC} $$
from which we see that the first vertical map is also an isomorphism.
This shows that $P$ is a generator. 
\smallskip

(7) By (6), $_A\cC$ is flat and $P$ is a generator.
Consider any epimorphism $f:M\to N$ in $\cM^\cC$. From this we obtain the
commutative diagram
$$\xymatrix{ 
  \Hom^\cC(P,M)\ot_SP \ar[d]^\simeq \ar[rr]^{\Hom(P,f)\ot\id_P} &  &
  \Hom^\cC(P,N)\ot_SP \ar[d]^\simeq\ar[r] & 0 \\
  M \ar[rr] & & N \ar[r] & 0, }$$
where the vertical maps are isomorphisms by \ref{gen} and hence 
the exactness of the bottom row implies exactness of the top row.
Now faithulness of the functor $-\ot_SP$ implies that $\Hom^\cC(P,f)$
is an epimorphism and hence $P$ is projective in $\cM^\cC$.
\endproof

\begin{thm}{\bf Remark.} \em Notice that the condition $_TP$ finitely
generated is satisfied if $P$ is a generator in $\cM_A$. 
For Galois comodules $P$ this is 
the case provided $\vareps:\cC\to A$ is surjective.
\end{thm}

\begin{thm}{\bf Corollary.}
Let $P\in \cM^\cC$ be a Galois comodule.
\begin{zlist}
\item If 
\begin{tabular}[t]{cl}
{\rm(i)} & $_SP$ is projective or \\
  {\rm(ii)}& $_SP$ is locally projective and $_TP$ is finitely generated, 
\end{tabular}\\
then $\cM^\cC$ is equivalent to the full category of $_{^*\cC}\cM$   
subgenerated by the ${^*\cC}$-module $P$, i.e., $\cM^\cC=\sigma[{_{^*\cC}P}]$.

\item If $_SP$ is finitely generated and projective, then  
      $\cM^\cC={_{^*\cC}\cM}$.
\end{zlist}
\end{thm}
\proof (1) Under the given conditions, $_A\cC$ is locally projective 
(see \ref{S.prop}(3),(4)) and 
  $\cM^\cC= \sigma[{_{^*\cC}\cC}]$. Since $P$
subgenerates $\cC$ it  is a subgenerator in $\cM^\cC$ and hence the assertion
follows. 

(2) The condition implies that $_A\cC$ is finitely generated and projective  
and hence $\cM^\cC={_{^*\cC}\cM}$ (by \cite[19.6]{BW}). 
\endproof 

\begin{thm}\Label{sem.gal}{\bf Semisimple Galois comodules.}
Assume $_A\cC$ to be flat. For a semisimple right $\cC$-comodule $P$,
the following are equivalent:
\begin{blist}
 \item $P$ is a Galois comodule;
 \item $P$ is a generator in $\cM^\cC$;
 \item $\mu_\cC:\Hom^\cC(P,\cC)\ot_SP\to \cC$ is surjective. 
\end{blist}
In this case $\cC$ is a right semisimple coring (and $_A\cC$ is projective).
\end{thm}
\proof Since $P$ is semisimple it is a generator in $\sP$ (see \ref{semis}).

(a)$\Ra$(c) This is trivial. 

(c)$\Ra$(b) Surjectivity of $\mu_\cC$ means that $\cC$ is $P$-generated.
 Since $\cC$ is a subgenerator in $\MC$ (see \cite[18.13(1)]{BW})
 this implies $\sP=\cM^\cC$.

(b)$\Ra$(a) follows from \ref{gen}(2).
\endproof 

\begin{thm}\Label{simp.gal}{\bf Simple Galois comodules.}
If $_A\cC$ is flat the following are equivalent:
\begin{blist}
 \item There is a simple Galois comodule in $\cM^\cC$;
 \item every non-zero comodule in $\cM^\cC$ is a Galois comodule;
 \item $\cC$ is homogenously semisimple as right comodule;
 \item $\cC$ is right semisimple and all simple right comodules are 
       isomorphic.  
\end{blist}
\end{thm}
\proof This follows by the characteriztions of simple right semisimple
corings in \cite[19.15]{BW} and the fact that each non-zero comodule is a
generator in this case. 
\endproof

\section{Galois comodules f.g. projective as $A$-modules} 

Some of the results in the preceding section were proved in
\cite[18.27]{BW} for the special case when $P_A$ is finitely generated and
 projective. As already observed (in \ref{pro.prop}) the latter condition 
provides nice properties of the  functor $\Hom^\cC(P,-)$ which 
 will lead to a left right symmetry of the Galois comodules.
If $P_A$ will be finitely generated and projective
 we denote a dual basis of $P$ by $p_1,\ldots,p_n \in P$ and
$\pi_1,\ldots,\pi_n \in P^*$.

\begin{thm}\Label{Cor.flat}{\bf $\mu_\cC$ splitting in $\MC$.}
Let $P\in \cM^\cC$ with $P_A$ finitely generated and projective and
$S=\End^\cC(P)$. Assume that
$$(P^* \ot_{S}P)\ot^\cC P^*\simeq P^*\ot_{S}(P\ot^\cC P^*)$$
canonically. Then the following are equivalent:
\begin{blist}
\item The map $\mu_C: \Hom^\cC(P,\cC)\ot_S P\to \cC$ is a splitting
epimorphism    in $\cM^\cC$; 
\item $\mu_\cC$ is an isomorphism. 
\end{blist}
The condition is satisfied provided $P$ is strongly $(\cC,A)$-injective,
or $P^*_S$ is flat, or $P^*$ is coflat or $(C,A)$-injective.   

\end{thm}
\proof We only have to prove (a)$\Ra$(b). 
It follows from \ref{pro.prop} that $P^*$ is
$P$-adstatic. By assumption, there is a splitting exact sequence in $\cM^\cC$,
 $$\xymatrix{
 0\ar[r]&  K \ar[r]& \Hom^\cC(P,\cC)\ot_S P \ar[r] & \cC \ar[r] &0 }.$$
Since  $P^*$ is $P$-adstsatic by \ref{pro.prop}(3), applying $\Hom^\cC(P,-)$
  yields an exact sequence 
$$\xymatrix{
0\ar[r]&\Hom^\cC(P,K)\ar[r]& \Hom^\cC(P,\Hom^\cC(P,\cC)\ot_SP)\ar[r]
^{\qquad\quad\simeq} & 
 \Hom^\cC(P,\cC) }.$$
From this we see $\Hom^\cC(P,K)=0$ and - since $K$ is a $P$-generated
comodule - this implies $K=0$.
\endproof

\begin{thm}\Label{left.right}{\bf Lemma.} Let $P\in \cM^\cC$ with 
$P_A$ finitely generated and projective. Then the following
are equivalent: 
\begin{blist}
\item $\cC$ is $P$-static as right $\cC$-comodule;
\item $\cC$ is $P^*$-static as left $\cC$-comodule.
\end{blist} 
\end{thm}
\proof 
The canonical map $\phi:P\to {^*(P^*)},\; \phi(p)(f)=f(p)$
for $p\in P$, $f\in P^*$, is bijective and the diagram 
$$\xymatrix{ 
   P^* \ot_S P\ar[d]^{\id\ot \phi}  \ar[r]& \cC \ar[d]^= &
       g \ot p \ar@{|->}[r] \ar@{|->}[d] &
              \sum g(p\sw{0}) p_\sw{1} \ar@{|->}[d] \\
 P^* \ot_S {^*(P^*)} \ar[r]& \cC & g \ot h \ar@{|->}[r] & 
  (\id_\cC\ot h)\varrho^{P^*}(g),
 } $$  
is commutative by the equalities 
 $$\begin{array}{c}
(\id_\cC\ot \phi(p))\varrho^{P^*}(g)= 
  \sum_i (g\ot \id_\cC)\varrho^P(p_i) \phi(p)(\pi_i)
  =  (g\ot \id_\cC)\varrho^P(\sum_i p_i \pi_i(p))=\sum g(p\sw{0}) p_\sw{1}.
\end{array}$$
By definition, $\cC$ is $P$-static provided the map in the top row is an
isomorphism of right $\cC$-comodules, and $\cC$ is $P^*$-static as left
$\cC$-comodule provided the map in the bottom row of the diagram is an
isomorphism of left $\cC$-comodules.   
\endproof 

Recall that for any bimodule $_BP_A$ with $P_A$  finitely generated and
projective (with dual basis as above), the 
 $(A,A)$-bimodule $P^*\ot_BP$ is an $A$-coring with   
coproduct and counit defined by   
 $$ \begin{array}{c}
 \uDelta:P^*\ot_BP \to (P^*\ot_B P)\ot_A (P^*\ot_BP), \quad 
      f\ot p \mapsto \sum f\ot p_i \ot \pi_i \ot p,\\[+2mm]
 \ueps:  P^*\ot_BP \to A,\quad f\ot p \mapsto f(p).
\end{array}$$

For a Galois comodule this coring is isomorphic to $\cC$.

\begin{thm}\Label{gal.Afg} {\bf Galois comodules with $P_A$ f.g. projective.} 
Let $P\in \cM^\cC$ with $P_A$ finitely generated and projective
and $S=\End^\cC(P)$. Then the following are equivalent:
\begin{blist}
\item $P$ is a Galois right $\cC$-comodule;
\item $\cC$ is $P$-static as right $\cC$-comodule;
\item $P^*$ is a Galois left $\cC$-comodule;
\item $\cC$ is $P^*$-static as left $\cC$-comodule;
\item $\tilde\mu_A: P^*\ot_S P \to \cC$ is an $A$-coring isomorphism.
\end{blist}
\end{thm}
\proof (a)$\LRa$(b) This is shown in \cite[18.26]{BW}.  
\smallskip

(c)$\LRa$(d) The assertion is the left hand version of (a)$\LRa$(b). 
\smallskip

(b)$\LRa$(d) This is proved in \ref{left.right}.
\smallskip 

(b)$\LRa$(e) It remains to show that $\tilde\mu_A$ is a coring morphism. Proofs
for this are given in 
\cite[Proposition 2.7]{ElGo.com} and \cite[18.26]{BW}.
With our notation it is seen by the following argument.
  For any $p\in P$ and $f\in P^*$, $p= \sum_i p_i \pi_i(p)$,  
 $$\begin{array}{rl}
 \tilde\mu_A(f\ot p) = \sum f(p\sw{0})p\sw{1} = 
   \sum  \sum_i f({p_i}\sw{0}){p_i}\sw{1} \pi_i(p),\;\mbox{ and }
 \end{array}$$
 $$\begin{array}{rl}
 (\tilde\mu_A\ot \tilde\mu_A)\circ \uDelta (f\ot p)& = 
  \sum \sum_i f({p_i}\sw{0}) {p_i}\sw{1} \ot \pi_i(p\sw{0}) p\sw{1} \\[+1mm]
   & = \sum f(p\sw{00})p\sw{01}\ot  p\sw{1} \\[+1mm]
   & = \sum f(p\sw{0})p\sw{11}\ot  p\sw{12} \\[+1mm]
   & = \Delta \circ \tilde\mu_A  (f\ot p) ,
 \end{array} $$
and it is easy to see that $\vareps\circ \tilde\mu_A = \ueps$. 
\endproof

\begin{thm}{\bf Remark.} \em 
It was shown in \ref{S.prop}(5) that for a Galois comodule $P\in \cM^\cC$,
$_SP$ locally projective and $_TP$ finitely generated, implies that 
$_A\cC$ is locally projective. In case $P_A$ is finitely generated and
projective, $_SP$ locally projective implies $_A\cC$ locally projective
without the additional assumption that $_TP$ is finitely generated
(see \cite[19.7]{BW}).  
\end{thm}
 
In the special situation that 
 $A$ is a $\cC$-comodule, i.e., there is a grouplike element 
$g\in\cC$,  and $S=\End^\cC(A)$, it is a Galois (right) comodule 
($(\cC,g)$ is a Galois coring) if and only if the map (compare introduction)
     $$A\ot_S A\to \cC,\; a\ot a'\mapsto aga',$$
is an isomorphism. Under the given conditions, $A\ot_S A$ has a 
canonical  coring structure (Sweedler coring) and the map is a coring
isomorphisms (see \cite[28.18]{BW}).   

 At various places we have observed nice properties of strongly $(\cC,A)$-injective 
comodules. For Galois comodules this notion is symmetric in the 
following sense - an observation also proved in \cite[Theorem 7.2]{Brz:gal}. 

\begin{thm}\Label{cov.symm} {\bf Strongly $(\cC,A)$-injective Galois
comodules.} 
Let $P$ be a Galois comodule with $P_A$ finitely generated and
projective and $S=\End^\cC(P)$. Then the following are equivalent:
\begin{blist}
\item $P$ is strongly $(\cC,A)$-injective;
\item $P^*$ is strongly $(\cC,A)$-injective;
\item the inclusion $i:S\to T$ is split by an $(S,S)$-bilinear map.
\end{blist}
\end{thm}
\proof This follows from \ref{CA-inj} and symmetry. 
\endproof

\begin{thm}\Label{P.stat}{\bf $P$-static comodules.} 
Let $P\in \MC$ with $P_A$ 
finitely generated and projective and assume $\cC$ to be $P$-static. 
Then 
 $N\in \cM^\cC$ is $P$-static, provided   
$$(N\ot^\cC P^*) \ot_S P \simeq  N\ot^\cC (P^* \ot_S P) $$
 canonically.
This holds if $N$ is $(\cC,A)$-injective, or $N\otC -$ is right exact,
or $P$ is strongly $(\cC,A)$-injective, or $P$ is  flat as $S$-module.
\end{thm}
\proof The first claim follows by the isomorphisms
$$\Hom^\cC(P,N)\ot_S P \simeq (N\ot^\cC P^*) \ot_S P \simeq N\ot^\cC (P^* \ot_S P)
 \simeq N.$$
The remaining assertions are derived from \ref{ass.cot}.
\endproof

The relevance of the isomorphism in \ref{P.stat} was also 
observed in \cite[Proposition 2.4]{CaGrVe}.
Notice that \ref{P.stat} shows again - in this special case - that
$P$ is a Galois comodule provided $\cC$ is
 $P$-static (see \ref{gal.Afg}), and that Galois comodules are
generators provided they are flat over their endomorphism rings (see
\ref{S.prop}). 

\begin{thm}\Label{equival}{\bf Equivalences.} Let $P\in \cM^\cC$ be a 
Galois comodule with $P_A$ finitely generated and projective.
Then
$$\Hom^\cC(P,-): \cM^\cC\to \cM_S$$ 
is an equivalence with inverse functor $-\ot_SP$ provided that
\begin{rlist}
 \item $P$ is strongly $(\cC,A)$-injective,  or  
 \item $P^*$ is $(\cC,A)$-injective and $_SP$ is flat, or 
 \item $P^*$ is coflat and $_SP$ is flat, or 
 \item $\cC$ is a coseparable coring. 
\end{rlist}
\end{thm}
\proof  
 Under each of the  conditions (i)-(iii) all right $S$-modules are
$P$-adstatic by \ref{pro.prop} and the right $\cC$-comodules are $P$-static by
\ref{P.stat}.

(iv) For a coseparable coring all comodules are strongly
 $(\cC,A)$-injective and hence (i) holds.  
\endproof

Parts of the preceding theorem are proved in \cite[Proposition 7.3]{Brz:gal}. 
Here we offer alternative proofs and do not require $_A\cC$ to be flat 
in the first case.

\begin{thm}\Label{rem.Ca}{\bf Remarks.} \em
  In \cite[Proposition 5.6]{CaGrVe}, the coring $\cC$ is required to be 
coseparable, finitely generated and projective as right $A$-module, 
and $\mu_\cC:\Hom^\cC(P,\cC)\ot_S P\to \cC$ should be surjective. 
These conditions immediately imply that $\mu_\cC$ splits in $\cM^\cC$ and
that $P^*$ is coflat. Hence  $P$ is a Galois comodule by \ref{Cor.flat} 
and the claim of \cite[Proposition 5.6]{CaGrVe} - namely that  $\Hom^\cC(P,-)$
is an equivalence - follows from \ref{equival}. 
\end{thm}

\begin{thm}\Label{flat.sub}{\bf Splitting over a subring of $\End^\cC(P)$.}
Let $P\in \cM^\cC$ with $P_A$ finitely generated and projective, 
$S=\End^\cC(P)$ and $B\subseteq S$ a subring. Assume $P^*$ to be flat as 
a right $B$-module, or $P$ to be $B$-strongly $(\cC,A)$-injective. Then the
following are equivalent: 
\begin{blist}
\item The canonical map $\mu_\cC':\Hom^\cC(P,\cC)\ot_B P \to \cC$ is a
     splitting epimorphism in $\cM^\cC$;
\item $P^*$ is a Galois comodule and is $(S,B)$-projective as right module. 
\end{blist}
\end{thm}
\proof
 (a)$\Ra$(b) Cotensoring with $-\ot^C P^*$, $\mu_\cC'$ yields the splitting
epimorphism in $\cM_S$,
  $$ P^*\ot_B S \simeq (P^*\ot_B P)\ot^\cC P^*\to \cC\ot^\cC P^*\simeq P^*,$$ 
where the first isomorphism is due to the conditions on $P^*$ or $P$
 (see \ref{hom.cot}).
This shows that $P^*$ is $(S,B)$-projective (see \cite[20.3]{WiBim}).  
In particular, $P^*$ is flat as $S$-module. 

Furthermore, $\mu_\cC'$ factors over a splitting epimorphism 
$\mu_C:\Hom^\cC(P,\cC)\to \cC$ in $\cM^\cC$. By Corollary \ref{Cor.flat}, this
implies that $\mu_\cC$ is an isomorphsm, i.e., $P$ is a Galois module.  
\smallskip

(b)$\Ra$(a)  By assmption, the map $P^*\ot_B S\to P^*$ splits in $\cM_S$.
Since $P^*$ (hence $P$) is a Galois comodule, 
tensoring with $-\ot_S P$ yields a splitting comodule epimorphism
    $$\begin{array}{c}
   P^*\ot_B P\simeq P^*\ot_B S\ot_S P \to P^*\ot_S P\simeq \cC. 
\end{array}$$ 
\endproof

\begin{thm}\Label{rem.Br}{\bf Remarks.} \em 
 (1) In case $A$ is an algebra over a
field (or a commutative von Neumann regular ring) $B$, 
then in \ref{flat.sub}, $P^*$ is always a flat $B$-module and the assertion  
yields  \cite[Theorem 4.4]{Brz:gal} as a special case.

(2) As pointed out in the introduction entwining structures can be 
considered as corings and hence the assertions in \ref{pro.prop}
and \ref{CA-inj} may be compared with 
Lemma 4.1 and Remarks 4.2 and 5.3 in \cite{SchSch}. Furthermore,
the splitting properties considered in \ref{flat.sub} are 
related to Remark 4.4, Theorem 2.2 and results of Section 5 
in \cite{SchSch}.  
\end{thm}

\section{Direct sums of f.g. projective $A$-modules.}

For the investigation of direct sums of modules the following technical 
observation is helpful.
For a direct sum of modules $P=\bigoplus_\Lambda P_\lambda$, denote by
$\epsilon_\lambda: P_\lambda\to P$ and $\pi_\lambda:P\to P_\lambda$
the canonical injections and projections. Recall that 
the identity of $P$ can be written as the formal sum 
$\sum_\Lambda \epsilon_\lambda\circ\pi_\lambda$.

\begin{thm}\Label{sum.lemma}{\bf Lemma.}
Let $P=\bigoplus_\Lambda P_\lambda$ be a direct sum of right $A$-modules
and $S\subseteq \End_A(P)$ a subring containing   
 $\epsilon_\lambda\circ\pi_\lambda$, for each $\lambda\in \Lambda$.
 Then, for any $K\in \cM_A$,
$$\begin{array}{c}
\Hom_A(P,K)\ot_SP \simeq (\bigoplus_\Lambda\Hom_A(P_\lambda,K))\ot_SP.
\end{array}$$  
\end{thm}
\proof  
Clearly the inclusion $(\bigoplus_\Lambda\Hom_A(P_\lambda,K))\ot_SP \to
\Hom_A(P,K)\ot_SP$ is injective. To see that it is surjective take any 
$f\in \Hom_A(P,K)$, $p\in P$, and write
$$\begin{array}{c}
 f\ot p = f\ot (\sum_\Lambda \epsilon_\lambda\circ\pi_\lambda(p))
         = \sum_\Lambda f\circ\epsilon_\lambda\circ\pi_\lambda \ot p, 
\end{array}$$
where  $f\circ\epsilon_\lambda\circ\pi_\lambda\in  \Hom_A(P,K)$
and  $f\circ\epsilon_\lambda \in  \Hom_A(P_\lambda,K)$.
\endproof
 
With this isomorphism the special structure of comodules that are finitely
generated as $A$-modules can be extended to direct sums of modules of this
type.  

\begin{thm}\Label{left.right.s}{\bf $P^*\ot_SP$ as left comodule.} \em
Consider a family $\{P_\lambda\}_\Lambda$ of comodules
$P_\lambda\in \cM^\cC$ such that each $P_\lambda$ is finitely generated and
projective as $A$-module. 
Then $P=\bigoplus_\Lambda P_\lambda$ is in $\cM^\cC$. 
Since all $P^*_\lambda$ are left $\cC$-comodules, 
their direct sum $\bigoplus_\Lambda P^*_\lambda$ is a left $\cC$-comodule. 
For $S=\End^\cC(P)$, \ref{sum.lemma} yields the identification 
  $P^*\ot_S P \simeq (\bigoplus_\Lambda P^*_\lambda)\ot_S P$
which makes $P^*\ot_S P$ to a left $\cC$-comodule. 

{\em
The following are equivalent: 
\begin{blist}
\item $\cC$ is $P$-static as right $\cC$-comodule;
\item $\cC$ is $\bigoplus_\Lambda P_\lambda^*$-static as left $\cC$-comodule.
\end{blist} 
}
\end{thm}
\proof Applying \ref{sum.lemma} repeatedly yields isomorphisms
 $$\begin{array}{rcl}
  P^*\ot_S P &\simeq &
   (\bigoplus_\Lambda P_\lambda^*)\ot_S (\bigoplus_\Lambda P_\lambda)
\\[+1mm]
     &\simeq & (\bigoplus_\Lambda P_\lambda^*)\ot_S (\bigoplus_\Lambda
{^*(P_\lambda^*)}) \\[+1mm]  
 &\simeq & (\bigoplus_\Lambda P_\lambda^*)\ot_S
{^*(\bigoplus_\Lambda(P_\lambda)^*)} .  
\end{array} $$
With this isomorphisms the proof of \ref{left.right} applies.
\endproof

The construction of corings for finitely generated projective $A$-modules
can also be extended to direct sums of modules of this type.  
 
\begin{thm}\Label{cor.dir}{\bf Coring structure on direct sums.} \em 
Consider a family $\{{_{S_\lambda}P}_\lambda\}_\Lambda$ 
of $(S_\lambda,A)$-bimod\-ules that are finitely
generated and projective as right $A$-modules with dual basis 
$p_{\lambda_i}\in P_\lambda$, $\pi_{\lambda_i}\in P^*_\lambda$.
For each $\lambda\in \Lambda$ we have corings with coproducts  
$$\begin{array}{c}
\uDelta_\lambda:P_\lambda^*\ot_{S_\lambda} P_\lambda \to 
P_\lambda^*\ot_{S_\lambda}P_\lambda \ot_A P_\lambda^*\ot_{S_\lambda}P_\lambda, 
\quad f_\lambda\ot p_\lambda 
\mapsto \sum_i f\ot p_{\lambda_i} \ot\pi_{\lambda_i} \ot p, 
\end{array} $$
and the evaluation as counit. 

Put $P=\bigoplus_\Lambda P_\lambda$ and $S_\Lambda=\bigoplus_\Lambda
S_\lambda$. Then $S_\Lambda$ is a ring without unit and $P$ is a left
$S_\Lambda$ module by componentwise multiplication. This means in particular
that we can identify 
 $$P_\lambda^*\ot_{S_\lambda} P_\lambda = P_\lambda^*\ot_{S_\Lambda} P, $$
and so, by the universal property of the coproduct, 
the $\uDelta_\lambda$ yield an $(A,A)$-bilinear map
$$ \begin{array}{c}
 \uDelta:  (\bigoplus_\Lambda P^*_\lambda) \ot_{S_\Lambda}P \to  
    P^*\ot_{S_\Lambda}P\ot_A P^*\ot_{S_\Lambda} P.
\end{array}$$
By \ref{sum.lemma}, we may assume 
$P^*\ot_{S_\Lambda}P=(\bigoplus_\Lambda P^*_\lambda)\ot_{S_\Lambda} P$ 
and thus $\uDelta$ defines an $A$-coring structure on $P^*\ot_{S_\Lambda}P$
with the evaluation as counit. By construction, $P^*\ot_{S_\Lambda}P$ 
is isomorphic to the coring coproduct 
$\bigoplus_\Lambda P_\lambda^*\ot_{S_\lambda}P_\lambda$.

For any subring $S\subseteq T=\End_A(P)$ which contains $S_\Lambda$,
there is a canonical epimorphism $\beta:P^*\ot_{S_\Lambda}P\to P^*\ot_{S}P$
and we have the  maps
$$\xymatrix{
    0\ar[r] & P^*\ot_{S_\Lambda}P \ar[r]^{\uDelta\qquad\quad} \ar[d]^\beta & 
    P^*\ot_{S_\Lambda}P \ot_A P^*\ot_{S_\Lambda}P \ar[d]^{\beta\ot\beta}\\
    & P^*\ot_{S}P \ar@{-->}[r]^{\uDelta'\qquad\quad} 
    &P^*\ot_{S}P\ot_A P^*\ot_{S}P, }
$$
where ${\uDelta'}$ exists provided $\Ke \beta \subseteq \Ke (\beta\ot
\beta)\circ \uDelta$. To show this
take any $f\ot_{S_\Lambda} p\in P^*\ot_{S_\Lambda}P$ such that 
$ \sum_i f\ot_S p_i \ot_A \pi_i \ot_S p = 0$. 
Then $0= \sum_i f\ot_S p_i\pi_i(p) = f\ot_S p$ 
proving that $f\ot_S p\in \Ke\beta$.
It is easy to see that a similar argument works for finite sums of 
elements of the form $f\ot p$. This shows that our condition on $\Ke \beta$ is
satisfied and that ${\uDelta'}$ exists making  $P^*\ot_{S}P$ an $A$-coring
with the evaluation map as counit.  

The coring structure on  $P^*\ot_SP$ as given here was introduced 
in \cite{ElGo} (along a different line of arguments) and these corings are
called {\em infinite comatrix corings} there. 
\end{thm}

With this preparation we are now able to extend the characterization 
of Galois comodules which are finitely generated and projective to  
those which are direct sums of such comodules.   

\begin{thm}\Label{gal.proj} {\bf $P_A$ direct sum of f.g. projectives.}
Consider a family $\{P_\lambda\}_\Lambda$ 
of $\cC$-comodules that are finitely
generated and projective as right $A$-modules (with dual basis
as in \ref{cor.dir}) and put $S_\lambda=\End^\cC(P_\lambda)$.
Then for $P=\bigoplus_\Lambda P_\lambda$ and $S=\End^\cC(P)$,
the following are equivalent: 

\begin{blist}
\item $P$ is a Galois right $\cC$-comodule;
\item $\cC$ is $P$-static as right $\cC$-comodule;
\item $\bigoplus_\Lambda P_\lambda^*$ is a Galois left $\cC$-comodule;
\item $\cC$ is $\bigoplus_\Lambda P_\lambda^*$-static as left $\cC$-comodule;
\item $\tilde\mu_A: P^*\ot_S P \to \cC$ is an $A$-coring isomorphism.
\end{blist}
\end{thm}
\proof (a)$\LRa$(b)  
One implication is trivial. Assume $\cC$ to be $P$-static. Then for any 
$K\in \cM_A$, the isomorphisms from \ref{sum.lemma} yield 
 $$\begin{array}{rcl}
 \Hom_A(P,K)\ot_SP &=& \bigoplus_\Lambda \Hom_A(P_\lambda,K)\ot_SP \\[+1mm]
    & \simeq& \bigoplus_\Lambda  K\ot_A P^*_\lambda \ot_SP \\[+1mm]
    & \simeq & K\ot_A (\bigoplus_\Lambda P^*_\lambda) \ot_SP \\[+1mm]
     &\simeq & K\ot_A P^*\ot_SP \simeq K\ot_A\cC.
\end{array}$$
Now the assertion follows from  \ref{gal.comod}(c).
\smallskip 

 (c)$\LRa$(d) is the left hand version of (a)$\LRa$(b).
\smallskip

(b)$\LRa$(d) This is shown in \ref{left.right.s}.
\smallskip

(b)$\LRa$(e) It remains to show that $\tilde\mu_A$ is a coring 
morphism. As mentioned in the proof of \ref{gal.Afg}, the maps
$P_\lambda^*\ot_{S_\lambda} P_\lambda\to \cC$ are coring morphisms. 
Hence $P^*\ot_{S_\Lambda} P\to \cC$ is a coring morphism and so is the 
factorisation $P^*\ot_{S} P\to \cC$. Notice that this is also proved
in \cite[Lemma 3.7]{ElGo}.
\endproof

For any module $M$ that is a direct sum of finitely generated modules it is
convenient to consider the functor $\wHom(M,-)$. This is, for example,
outlined  in \cite[Section 51]{WiBu} and a straightforward transfer of the
related notions to comodules yields the following.     

\begin{thm}\Label{hat.functor} {\bf The functor $\wHomC(P,-)$.} \em
Given $P=\bigoplus_\Lambda P_\lambda$ as a direct sum of comodules that are  
finitely generated as $A$-modules, consider the morphisms for $N\in \MC$,
$$\wHomC(P,N) = \{f\in \HomC(P,N)\, |\, f(P_\lambda)=0\;
 \mbox{ for almost all } \lambda\in \Lambda\}.$$  
Then $\wS=\wHomC(P,P)$ is a subring - in fact a left ideal - in $S=\EndC(P)$   
with enough idempotents.  This induces a functor
$$\wHomC(P,-):\MC\to \cM_{\wS}$$ 
where $\cM_{\wS}$ is the category of
all right $\wS$-modules $X$ with $X\wS=X$.
There is a functorial isomorphism 
$\wHomC(P,-)\simeq \Hom^\cC(P,-)\ot_S \wS$ yielding  
the isomorphisms
 $$\wHomC(P,N)\ot_\wS P \simeq \Hom^\cC(P,N)\ot_S \wS \ot_\wS P  
  \simeq \Hom^\cC(P,N)\ot_S P,$$
and hence $N$ is $P$-static if and only if
 $\wHomC(P,N)\ot_\wS P\simeq N$.

If the $P_\lambda$ are finitely generated and projective as $A$-modules, then
 $$\begin{array}{c}
 \wHomC(P,-) \simeq \bigoplus_\Lambda \HomC(P_\lambda,-)\simeq
      -\otC (\bigoplus_\Lambda P_\lambda^*).
\end{array}$$
So in this case the functor  $-\otC (\bigoplus_\Lambda P_\lambda^*)$ is right 
adjoint to the functor $-\ot_S P$ (thus yielding 
\cite[Proposition 4.5]{ElGo}). 
\end{thm}

\begin{thm}\Label{P.stat.s}{\bf $P$-static comodules.} 
Let $P=\bigoplus_\Lambda P_\lambda$ where the $P_\lambda\in \MC$ 
are finitely generated and projective as $A$-modules and 
assume $\cC$ to be $P$-static. Then 
 $N\in \cM^\cC$ is $P$-static, provided   
$$\begin{array}{c}
(N\ot^\cC \bigoplus_\Lambda P_\lambda^*) \ot_S P \simeq  
N\ot^\cC (\bigoplus_\Lambda P_\lambda^* \ot_S P) 
\end{array}$$  
canonically.
This holds if $N$ is $(\cC,A)$-injective, or $N\otC -$ is right exact,
or the $P_\lambda$'s are strongly $(\cC,A)$-injective, or $P$ is flat as
$S$-module. 
\end{thm}
\proof The first claim follows by the isomorphisms
$$\begin{array}{rcl}
 \wHomC(P,N)\ot_\wS P 
 &\simeq& (\bigoplus_\Lambda \HomC(P_\lambda,N)) \ot_S P \\[+2mm]
 &\simeq&  \bigoplus_\Lambda (N\ot^\cC P_\lambda^*) \ot_S P \\[+2mm]
 &\simeq&  N\ot^\cC (\bigoplus_\Lambda P^*_\lambda \ot_S P) \, \simeq\, N.
\end{array}$$
The remaining assertions are derived from \ref{ass.cot}.
\endproof

 \begin{thm}\Label{equival.P}{\bf Equivalences.} 
Let $P=\bigoplus_\Lambda P_\lambda$ be a Galois comodule, where the
$P_\lambda\in \MC$ are finitely generated and projective as $A$-modules.
Then 
$$\wHomC(P,-): \cM^\cC\to \cM_\wS$$ 
is an equivalence with inverse functor $-\ot_\wS P$ provided that
\begin{rlist}
 \item $P$ is strongly $(\cC,A)$-injective,  or  
 \item each $P_\lambda^*$ is $(\cC,A)$-injective and $_SP$ is flat, or 
 \item each $P_\lambda^*$ is coflat and $_SP$ is flat, or 
 \item $\cC$ is a coseparable coring. 
\end{rlist}
\end{thm}
\proof  The same arguments as for the proof of \ref{equival} apply.   
\endproof

By the isomorphisms $\HomC(P,-)\simeq -\otC P^*$, $P^*$ is coflat 
if and only if  $P$ is projective in $\MC$. Hence, by \ref{S.prop}, 
the condition (iii) in \ref{equival.P} implies that $_A\cC$ is flat
and that $P$ is a projective generator in the Grothendieck category $\MC$.
 With familiar arguments from module theory (see \cite[51.11]{WiBu})
this situation can be described in the following way.

\begin{thm}\Label{proj.gen} {\bf Projective generators in $\MC$.}
Let $P=\bigoplus_\Lambda P_\lambda$, where the
$P_\lambda\in \MC$ are finitely generated and projective as $A$-modules.
Then the following are equivalent:
\begin{blist}
\item $_A\cC$ is flat and $P$ is a projective generator in $\MC$;
\item $P$ is a Galois comodule and $_SP$ is faithfully flat;
\item $_A\cC$ is flat and $\wHomC(P,-): \cM^\cC\to \cM_\wS$ is an equivalence;
\item $_A\cC$ is flat and $-\ot_\wS P: \cM_\wS \to\cM^\cC$ is an
      equivalence.
\end{blist}
\end{thm}

Notice that similar characterisations are also proved in \cite[Theorem
4.7]{ElGo}. \medskip

Finally we ask when $\cC$ is a right Galois comodule. 
For this recall that $\End^\cC(\cC)\simeq \cC^*$ and that - in our notation -
 $\cC^*$  acts on $\cC$ from the right. Then the evaluation map 
  $$\mu_\cC: \cC \ot_{\cC^*}\Hom^\cC(\cC,\cC)\to \cC$$  
is an isomorphism and hence we conclude from \ref{prop.alpha}, \ref{gal.Afg}
and \ref{gal.proj}: 

\begin{thm}\Label{C.gal} {\bf $\cC$ as Galois comodule.}
  
\begin{zlist} 
 \item If $\cC_A$ is finitely generated and projective, 
    then $\cC$ is a Galois right $\cC$-comodule, $\cC^*$ is a Galois left
    $\cC$-comodule, and $K\ot_A\cC\simeq \Hom_A({^*\cC},K)$, 
    for any $K\in\cM_A$.

\item If $\cC=\bigoplus_\Lambda \cC_\lambda$ with right
       subcomodules $\cC_\lambda$ that are finitely generated and
       projective as right $A$-modules, then $\cC$ is a Galois right
       $\cC$-comodule and  $\bigoplus_\Lambda \cC^*_\lambda$
        is a Galois left $\cC$-comodule. 
\end{zlist}
\end{thm}

Notice that in (2) $\cC^*$ is not a left $\cC$-comodule unless the 
sum is finite (finiteness theorem, \cite[19.12]{BW}).
For further properties of $\cC$ as Galois comodule we refer to 
\cite[Section 7]{CaGrVe}.
\bigskip

{\bf Acknowledgement.} The author is very grateful to Tomasz 
Brzezi\'nski for an inspiring exchange of ideas and helpful comments.

\noindent 
{\bf Address:}\\
 Department of Mathematics \\
 Heinrich Heine University\\ 
 40225 D\"{u}sseldorf, Germany \\[+1mm]
 e-mail: wisbauer@math.uni-duesseldorf.de
              

\begin{thebibliography}{99}

 \bibitem{Brz:mod} Brzezi\'nski, T., {\em On modules associated to
        coalgebra-Galois extensions}, 
        J.\ Algebra 215, 290--317 (1999)
\bibitem{Brz:str} Brzezi\'nski, T., {\em The structure of corings. 
    	Induction functors, Maschke-type theorem, and Frobenius and 
    	Galois-type properties},    Alg. Rep. Theory 5, 389--410 (2002)
\bibitem{Brz:gal} Brzezi\'nski, T., {\em Galois Comodules},
   arXiv:math.RA/0312159v3 (2004)

 \bibitem{BW} Brzezi\'nski, T., Wisbauer, R.,
    {\em Corings and comodules},
  London Math. Soc. LNS 309, Cambridge University Press (2003)
 
\bibitem{CaGrVe}  Caenepeel, S., De Groot. E.,  Vercruysse, J., 
     {\em Galois theory for comatrix corings: Descent theory,
      Morita theory, Frobenius and separability properties}, \\
      arXiv:math.RA/0406436 (2004)

%\bibitem{Caene}  Caenepeel, S., 
%  {\em Galois corings from the descent theory point of view}, \\
% arXiv:math.RA/0311377 (2003)

%\bibitem{CaeIon:fro}  Caenepeel, S.,  Ion,  B., Militaru, G.,
% {\em The structure of Frobenius algebras and separable algebras},  
% K-Theory 19, 365--402 (2000)

\bibitem{ElGo.com} El Kaoutit, L., G\'omez-Torrecillas, J.,
   {\em Comatrix corings: Galois corings, descent theory, and a structure 
   theorem for cosemisimple corings}, Math. Z. 244, 887--906 (2003)

\bibitem{ElGo} El Kaoutit, L., G\'omez-Torrecillas, J., {\em Infinite Comatrix
   Corings}, \\
 arXiv:math.RA/0403249v1 (2004) 

\bibitem{SchSch} Schauenburg, P., Schneider, H.-J.,
   {\em On generalized Hopf Galois extensions}, \\
   arXiv:math.QA/0405184 (2004)

\bibitem{WiBu}  Wisbauer, R., 
      {\it Foundations of Module and Ring Theory},
      Gordon and Breach, Reading (1991) 

\bibitem{WiBim}  Wisbauer, R., 
    {\em Modules and Algebras: Bimodule Structure and Group Actions on 
    Algebras}, 
    Pitman Mono. PAM 81, Addison Wesley, Longman, Essex (1996)

\bibitem{WiTil} Wisbauer, R.,  {\em Tilting in module categories},
    in Abelian Groups, Module theory, and Topology,  Dikranjan and 
   Salce (Eds.),  Marcel Dekker LNPAM 201, 421--444 (1998)

\bibitem{WiStat} Wisbauer, R.,  {\em Static modules and equivalences}, 
      in Interactions Between Ring Theory and Representations of
         Algebras, van Oystaeyen and Saorin  (Eds.), 
    Marcel Dekker, 423--449 (2000)  

\bibitem{WiGal} Wisbauer, R., {\em On Galois Corings},  Hopf algebras 
     in non-commutative geometry and physics, S. Caenepeel and F. van
     Oystaeyen (eds.), LNPAM 239, Marcel Dekker, New York (2004) 

\bibitem{Zim} Zimmermann, W., {\em Modules with chain conditions for 
     finite matrix subgroups}, 
     J. Algebra 190, 68--87 (1997)    
 \end{thebibliography}
\end{document}